\documentclass[11pt]{amsart}
\usepackage{fullpage}
\usepackage{amsfonts,amssymb,mathrsfs}
\usepackage{graphics,color,hyperref}
\usepackage{paralist}
\usepackage{tikz}
\usepackage{barycentric}
\usetikzlibrary{arrows}
\usetikzlibrary{automata}

\DeclareMathAlphabet{\mathbbold}{U}{bbold}{m}{n}
\newtheorem{theorem}{Theorem}
\newtheorem*{theoremalph}{Theorem}
\newtheorem{conj}{Conjecture}
\newtheorem{prop}{Proposition}
\newtheorem{lemma}{Lemma}
\newtheorem{coro}{Corollary}
\theoremstyle{definition}
\newtheorem{definition}{Definition}
\theoremstyle{remark}
\newtheorem{example}{Example}
\newtheorem{remark}{Remark}
\newcommand{\zero}{\mathbbold{0}}
\newcommand{\unit}{\mathbbold{1}}
\def\rmax{\mathbb{R}_{\max}}
\newcommand{\smax}{\mathbb{S}_{\max}}
\newcommand{\sgn}{\operatorname{sgn}}
\newcommand{\R}{\mathbb{R}}
\newcommand{\N}{\mathbb{N}}

\newcommand{\sU}{\mathcal{U}}
\newcommand{\sC}{\mathcal{C}}
\newcommand{\sK}{\mathcal{K}}
\newcommand{\sT}{\mathcal{T}}
\newcommand{\floor}[1]{\lfloor #1\rfloor}

\newcommand{\mychoose}[2]{{{#1}\choose{#2}}}
\newcommand{\rbar}{\R\cup\{-\infty\}}
\newcommand\ntrop{N^{\text{tpath}}}
\newcommand\ntropext{N^{\text{trop}}}
\newcommand\nclass{N^{\text{path}}}
\newcommand\adj{\rm{adj}}
\newcommand\npm{N^{+-}}
\newcommand\nmp{N^{-+}}
\newcommand\pmp{-+}
\newcommand\ppm{+-}
\newcommand{\argmax}{\arg\max}
\newcommand{\sP}{\mathscr{P}}
\newcommand{\mon}{-}
\newcommand{\mop}{+}

\title{The number of extreme points of tropical polyhedra}
\author{Xavier {A}llamigeon}
\address{EADS Innovation Works, SE/IA -- Suresnes, France and CEA, LIST MeASI -- Gif-sur-Yvette, France}
\email{xavier.allamigeon@eads.net}
\author{{S}t{\'e}phane {G}aubert}
\address{INRIA and CMAP, \'Ecole Polytechnique, 91128 Palaiseau Cedex France}
\email{stephane.gaubert@inria.fr}
\author{Ricardo D. Katz}
\address{CONICET. Postal address:\!\! Instituto de Matem\'atica 
``Beppo Levi'',\! Universidad Nacional de Rosario, 
Avenida Pellegrini 250, 2000 Rosario, Argentina.}
\email{rkatz@fceia.unr.edu.ar}
\date{June 18, 2009}
\begin{document}
\begin{abstract}
The celebrated upper bound theorem of McMullen
determines
the maximal number of extreme points of a polyhedron
in terms of its dimension and the number of constraints which define it,
showing that the maximum is attained by the polar of the cyclic polytope.
We show that the same bound is valid
in the tropical setting, up to a trivial modification.
Then, we study the natural candidates
to be the maximizing polyhedra, which are the polars of a family of cyclic polytopes equipped with a sign pattern. 
We construct bijections between the extreme points of these polars
and lattice paths depending on the sign pattern,
from which we deduce explicit
bounds for the number of extreme points, showing in particular that the upper
bound is asymptotically tight as the dimension tends to infinity, keeping
the number of constraints fixed. 
When transposed to the classical case, the previous constructions yield some lattice path generalizations of Gale's evenness criterion.
\end{abstract}
\maketitle
\renewcommand{\thefootnote}{}\footnotetext{The first two authors were partially supported by the Arpege programme of the French National Agency of Research (ANR), project ``ASOPT'', number ANR-08-SEGI-005. The second author was also partially supported by the Digiteo project DIM08 ``PASO'' number 3389.}
\renewcommand{\thefootnote}{\arabic{footnote}}

\section{Introduction}
A fundamental result in discrete convex geometry is McMullen's upper bound
theorem, which settled a conjecture of Motzkin. We restate it for completeness.
\begin{theoremalph}[{\cite{mcmullen70}}]
Among all polytopes in $\mathbb{R}^d$ with $p$ extreme points, the cyclic
polytope maximizes the number of faces of each dimension.
\end{theoremalph}
The reader is referred to~\cite{ziegler98,matousek} for more information.
Recall  that a cyclic polytope is the convex hull of $p$ 
distinct points on the moment curve $\{(t,t^2,\ldots,t^d)\mid t\in \mathbb{R}\}$.

In particular, the number of facets (faces of dimension $d-1$) is known
to be at most
\[
U(p,d):=\mychoose{p-\floor{d/2}}{\floor{d/2}}+ 
\mychoose{p-\floor{d/2}-1}{\floor{d/2}-1}
\qquad 
\text{ for $d$ even, and}
\]
\[
U(p,d):= 2\mychoose{p-\floor{d/2}-1}{\floor{d/2}}  
\qquad 
\text{ for $d$ odd.}
\]
By duality, the same upper bound applies to the number of extreme
points of a $d$-dimensional polytope defined as the intersection
of $p$ half-spaces.

In max-plus or tropical convexity, the addition and multiplication
are replaced by the maximum and the addition, respectively. This
unusual model of discrete convexity has been studied
under different names by several authors. 
These include K. Zimmermann~\cite{zimmerman77},
Cohen, Gaubert and Quadrat~\cite{cgq00,cgq02}, with
motivations from discrete event systems and optimal control~\cite{maxplus97,ccggq99}, Kolokoltsov, Litvinov, Maslov and Shpiz~\cite{maslovkolokoltsov95,litvinov00}, with motivations from variations calculus and quasi-classics asymptotics.
The field attracted a new attention after the work of Develin and Sturmfels~\cite{DS}, who connected it with current developments of tropical geometry,
and also showed an unexpected relation with the study of tree metrics
in phylogenetic analysis.  
This has been the source of a number of works of the same authors
and of Joswig and Yu, see in particular~\cite{joswig04,DevelinYu,JSY07}. 
Tropical convexity can also be studied from the general perspective of abstract
convexity, a point of view adopted by Singer, see~\cite{cgqs04,NiticaSinger07a} and also by Briec and Horvath~\cite{BriecHorvath04}. 
Some further works developing or applying tropical convexity 
include~\cite{blockyu06,GK06a,katz05,gauser,AGG08,GM08,invariant09}. 

The notion of extreme point carries over to the tropical setting~\cite{BSS,GK},
and so, we may ask whether McMullen's theorem,
or rather, its dual, concerning the number of extreme points,
admits a tropical analogue.

Our first result, which we establish in Section~\ref{sec-upperbound},
shows that a McMullen type bound is still valid in the tropical setting. 

\begin{theorem}\label{th-upperbound}
The number of extreme rays of a tropical cone in $(\rbar)^d$
defined as the intersection of $p$ tropical half-spaces
cannot exceed $U(p+d,d-1)$.
\end{theorem}
The number $p+d$ instead of $p$ for the number of constraints
can be explained intuitively: in loose terms,
in the tropical world, all the numbers
are ``positive'', so the bound is the same as for a polyhedral
cone of the same dimension in which $d$ positivity
constraints would have been added to the $p$ explicit ones.
The number $d-1$ instead of $d$ for the dimension reflects
the fact that we are dealing with cones, rather than with convex
sets.

The most natural idea of proof would be to tropicalize
the classical method, which relies on the 
$f$-vector theory. However, some pathological features of the notions of faces
of tropical polytopes make somehow uneasy the development
of a tropical analogue of this theory (see~\cite{DevelinYu} for a discussion on faces). 
So, we choose a different approach, and establish Theorem~\ref{th-upperbound}
as a corollary of the classical upper bound theorem, using a deformation
argument in which the tropical polyhedron is seen as a degenerate limit
of a sequence of classical polyhedra. 

In the classical case, the polar of a cyclic polytope with $p$ extreme
points maximizes the number of extreme points among all the polytopes
of dimension $d$ defined by $p$ inequalities. In the tropical
case, the notion of polar can be defined as well~\cite{katz08}. This leads
us to define a family of tropical generalizations of the cyclic
polytopes, in which a sign pattern is incorporated.  Our second
result is the following.

\begin{theorem}\label{th-allowed}
The extreme rays of the polar of a signed cyclic polyhedral cone are in one to one correspondence with
tropically allowed lattice paths.
\end{theorem}
The definition of tropically allowed lattice paths is given in
Section~\ref{SectionSigned}, in which this theorem is proved.
We also give a characterization of the extreme rays of the
classical (non-tropical) analogue of this polar (Theorem~\ref{th-class}),
showing
that there are fewer extreme rays
in the tropical case. The latter lattice path characterization
is intimately related with Gale's evenness criterion, as
shown in Theorem~\ref{theo-classical}.

Recall that in the classical case, 
a point of a polyhedron is extreme if and only if the gradients
of the constraints that it saturates form a family of full rank.
The comparison between Theorem~\ref{th-allowed} and Theorem~\ref{th-class}
reflects the fact that the same is not true in the tropical setting.
Indeed, the proof of Theorem~\ref{th-allowed}, which relies
on a Cramer type result due to M. Plus
~\cite{plus} (Theorem~\ref{th-cramer} below, see~\cite{AGG08b}
for a recent account and also~\cite{RGST} for an alternative approach
due to Richter-Gebert, Sturmfels, and Theobald) 
and on the characterization of extreme points obtained by Allamigeon, Gaubert and Goubault in~\cite{AGG09}
(see Theorem~\ref{CharactExtreme} below) shows that in the tropical case,
some additional minimality condition must be added to the classical
rank condition. This fundamental discrepancy explains why
there are fewer extreme points in the tropical case.

The analogy with the classical case
suggests the conjecture that
the maximal number of extreme rays of a tropical polyhedron 
defined by $p$ inequalities in dimension $d$ is attained
by the polar of a signed cyclic polyhedral cone. 
However, the only evidence for this we have to offer is
Theorem~\ref{th-general}, which shows that as in the classical case, the system
of maximizing constraints can be chosen to be in general position.
Then, the signed cyclic polyhedral cones somehow provide
the simplest models for such systems of constraints.

These considerations lead us to estimate the maximal number of generators
of the polar of a signed cyclic polyhedral cone,
for which we provide explicit lower and upper bounds in Section~\ref{sec-estimates}.

Finally, we note that Theorem~\ref{th-allowed} may seem surprising
(and perhaps even disappointing)
in the light of the developments of enumerative tropical geometry, following
the work of Mikhalkin~\cite{mikhalkin}.
A deep result there (Mikhalkin's correspondence theorem)
is that certain classical enumerative invariants
(the number of algebraic curves satisfying appropriate constraints)
can be computed from their tropical analogues, taking into account certain
multiplicities. The results of the present paper are limited to
the linear case, but concern inequalities instead of equalities. Theorem~\ref{th-allowed} shows that the most natural enumerative object concerning
inequalities, the number of extreme points, does not tropicalize. More 
precisely, its proof shows that when deforming a classical polyhedron
to obtain a tropical polyhedron, some of the classical extreme points 
degenerate in points which are no longer extreme in the tropical sense.

\section{Bounding the number of extreme points of a tropical polyhedron}\label{sec-upperbound}

The symbol $\rmax$ will denote the max-plus semiring,
which is the set $\rbar$ equipped with the addition $(a,b)\mapsto a\oplus b:=\max(a,b)$
and the multiplication $(a,b)\mapsto a b:=a+b$. The zero and unit elements
will be denoted by $\zero$ and $\unit$, respectively,
so $\zero=-\infty$ and $\unit =0$. If $A=(a_{ij})$
is a $p\times d$ matrix with entries in $\rmax$, the matrix vector
product $Ax$ is naturally defined for $x\in \rmax^d$, $(Ax)_i:=\bigoplus_{1\leq j\leq d} a_{ij}x_j$, which can be rewritten as $\max_{1\leq j\leq d} a_{ij}+x_j$
with the classical notation. 

In the present section, we will apply some asymptotic arguments, mixing
classical and max-plus algebra, and so we will mainly use the classical
notation. However, in the next section, we shall make an
intensive use of the ``max-plus'' notation, which will make clearer
some analogies with the classical case. 
In all cases, the reader will easily avoid any ambiguity from the context.

A subset $\sC$ of $\rmax^d$ is a {\em tropical (convex) cone} if 
\[
u,v\in \sC,\;\lambda,\mu\in \rmax\implies \lambda u \oplus \mu v\in \sC \enspace .
\]
Here, we denote by $\oplus$ the tropical sum of vectors, which is nothing but
the entrywise max, and we denote by $\lambda u$ the vector obtained
by multiplying in the tropical sense (i.e., adding) the scalar $\lambda$
by each entry of the vector $u$.

We say that a non-zero vector $u\in \sC$ is an {\em extreme generator} of $\sC$
if $u=v\oplus w$ with $v,w\in \sC$ implies $u=v$ or $u=w$. The set of scalar 
multiples of an extreme generator of $\sC$ is an {\em extreme ray} of $\sC$. 
A subset $\sU$ of a tropical cone $\sC$ is said to be a {\em generating family} 
of $\sC$ if any vector $x \in \sC$ can be expressed as 
$x=\oplus_{1\leq k\leq K} \lambda_k u_k$ for some $K\in \N$, 
where $\lambda_k \in \rmax $ and $u_k \in \sU$ for all $1\leq k\leq K$. 
A tropical cone is {\em finitely generated} if it has a finite generating 
family.

Let us recall the following tropical analogue of the Minkowski theorem,
established by  Gaubert and Katz~\cite{GK06a,GK} and
Butkovi\v{c}, Schneider and Sergeev~\cite{BSS}.
\begin{theorem}[Tropical Minkowski Theorem~\cite{BSS,GK}, 
see also~\cite{GK06a}]\label{TheoMinkowski}
A closed tropical cone is generated by its extreme rays. 
\end{theorem}
This applies in particular to finitely generated 
tropical cones, which are always closed~\cite{GK}.
Then, we get a refinement (with the added characterization
in terms of extreme rays) 
of an observation made by several authors including 
Moller~\cite{Mol}, Wagneur~\cite{Wag},
and Develin and Sturmfels~\cite{DS}, showing that a finitely generated
tropical cone has 
a ``basis'' (generating family with minimal cardinality)
which is unique up to the multiplication of its vectors
by possibly different scalars. 
Observe that every generating family of $\sC$ must contain at least one 
vector in each extreme ray of $\sC$.   

To establish Theorem~\ref{th-upperbound}, we shall think of tropical
convex cones as limits of classical convex cones along
an exponential deformation. Let $\beta>0$ denote a parameter,
and let $E_\beta$ denote
the map from $\rmax^d$ to $\R^d$ which sends
the vector $x=(x_j)$ to the vector $(\exp(\beta x_j))$.
We denote by $L_\beta$ the inverse map of $E_\beta$. 

We shall use repeatedly the following inequalities,
which hold for any vector $v\in \rmax^d$,

\begin{align}\label{ineq-main}
\max_{1\leq j\leq d} v_j \leq \beta^{-1}\log \big(\sum_{1\leq j\leq d} 
\exp(\beta v_j)\big)\leq \beta^{-1}\log d + \max_{1\leq j\leq d} v_j 
\enspace .
\end{align}

We now prove Theorem~\ref{th-upperbound}.
Consider the tropical cone $\sC$ of $\rmax^d$ defined as the intersection
of the following $p$ tropical half-spaces:
\begin{align}\label{e-def-C}
\max_{1\leq j\leq d} a_{ij}+x_j \leq \max_{1\leq j\leq d} b_{ij}+x_j\enspace ,\qquad
1\leq i\leq p \enspace ,
\end{align}
and let $\sC(\beta)\subset \R^d$ denote the (ordinary) convex cone
consisting of the vectors $y$ satisfying the inequalities:
\[
y_i \geq 0\enspace , \qquad  1\leq i\leq d \enspace ,
\]
\begin{align}\label{ineq-def}
\frac{1}{d}\sum_{1\leq j\leq d}\exp(\beta a_{ij}) y_j \leq \sum_{1\leq j\leq d} \exp(\beta b_{ij})y_j
\enspace ,\qquad 1\leq i\leq p \enspace .
\end{align}
If $x\in \sC$, then 
\begin{align*}
\frac{1}{d}\sum_{1\leq j\leq d}\exp(\beta (a_{ij}+x_j))
&\leq \exp(\beta(\max_{1\leq j\leq d} a_{ij}+x_j))\\
&\leq \exp(\beta(\max_{1\leq j\leq d} b_{ij}+x_j))\\
&\leq \sum_{1\leq j\leq d} \exp(\beta (b_{ij}+x_j)) \enspace, 
\end{align*}
which shows that $y:=E_\beta(x)$ belongs to $\sC(\beta)$.

Consider now the simplex
\[
\Sigma:=\{y\in \R^d\mid y\geq 0,\;\; \sum_{1\leq j\leq d} y_j =1\} \enspace .
\]
The extreme rays of the cone $\sC(\beta)$ are in one to one correspondence
with the extreme points of the convex set $\sC(\beta)\cap \Sigma$. By eliminating
the variable $y_d$, we identify the latter set with a convex subset
of $\R^{d-1}$ defined by $p+d$ affine inequalities. It follows
that the number $K(\beta)$ of extreme points of $\sC(\beta)\cap \Sigma$ is such that
\[
K(\beta) \leq U(p+d,d-1) \enspace .
\]
Let $\{u_k(\beta)\}_{k=1,\ldots ,K(\beta)}\subset \R^d$ 
denote a family obtained by ordering the extreme
points of $\sC(\beta)\cap \Sigma$ in an arbitrary way.

Since $u_k(\beta)\geq 0$, we can find a vector
$v_k(\beta) \in \rmax^d$ such that 
$u_k(\beta)=E_\beta(v_k(\beta))$.

Let us now fix a sequence $\beta_m$ tending to infinity.
Since $K(\beta)$ only takes a finite number of values, after 
replacing $\beta_m$ by a subsequence, we may assume
that $K:=K(\beta_m)$ is independent of $m$.

Let us consider an arbitrary index $k$ among $1,\ldots,K$.
Since $\sum_j (u_k(\beta))_j =1$ and $(u_k(\beta))_j\geq 0$, we deduce
that $\exp(\beta (v_k(\beta))_j)\leq 1$, and so, $v_k(\beta)$ belongs to the
set $[-\infty,0]^d$.
Since this set is compact, possibly after extracting $K$ subsequences
we may assume that for every index $1\leq k\leq K$,
$v_k(\beta_m)$ tends to some vector $v_k\in [-\infty,0]^d$ 
as $m$ tends to infinity. 

By applying the map $L_\beta$ to the relation 
$\sum_j \exp(\beta (v_k(\beta))_j)= \sum_j (u_k(\beta))_j =1$,
we get thanks to Inequality~\eqref{ineq-main}, 
\[
\max_{1\leq j\leq d} (v_k(\beta))_j \leq 0 \leq \beta^{-1}\log d + \max_{1\leq j\leq d} (v_k(\beta))_j \enspace, 
\]
and so
\begin{align}\label{EqualZero}
\max_{1\leq j\leq d} (v_k)_j =0 \enspace .
\end{align}

We claim that the family $\{v_k \}_{k=1,\dots ,K}$ 
generates the tropical cone $\sC$.

By setting $y=E_\beta(v_k(\beta))$ in Inequality~\eqref{ineq-def},
applying the order preserving map $L_\beta$ to both sides
of this expression,
and using Inequality~\eqref{ineq-main},
we get
\[
-\beta^{-1}\log d + \max_{1\leq j\leq d} a_{ij}+(v_k(\beta))_j \leq \beta^{-1}\log d
+ \max_{1\leq j\leq d} b_{ij} + (v_k(\beta))_j \enspace .
\]
Taking $\beta:=\beta_m$ and letting $m$ tend to infinity, we deduce
that 
\[
\max_{1\leq j\leq d} a_{ij}+ (v_k)_j \leq \max_{1\leq j\leq d} b_{ij}+(v_k)_j \enspace ,
\]
which shows that $v_k\in \sC$.

Consider now an arbitrary element $x\in \sC$. Since $u_k(\beta)$ generates
the convex cone $\sC(\beta)$, we can express
the vector $E_\beta(x)\in \sC(\beta)$ as a linear combination 
\begin{align}\label{e-fond}
E_\beta(x) =\sum_{1\leq k\leq K} \delta_k u_k(\beta) = 
\sum_{1\leq k\leq K} \delta_k E_\beta(v_k(\beta)) \enspace,
\end{align}
for some scalars $\delta_k\geq 0$, which can be written
as $\delta_k=\exp(\beta \lambda_k(\beta))$
for some $\lambda_k(\beta)\in \rmax$.

We deduce from~\eqref{e-fond} that
\[ E_\beta(x)\geq \delta_k E_\beta(v_k(\beta)) \enspace,
\]
and so, for all $j$,
\[
x_j \geq \lambda_k(\beta) + (v_k(\beta))_j \enspace .
\]
Choosing any index $j$ such that $(v_k)_j=0$, 
which exists by~\eqref{EqualZero}, 
we deduce that $\lambda_k(\beta_m)$ is bounded
from above as $m$ tends to infinity. Hence, after extracting a new
subsequence, we may assume that $\lambda_k(\beta_m)$ converges to
some scalar $\lambda_k\in \rmax$. Then, letting $\beta=\beta_m$
tend to infinity in~\eqref{e-fond}, and using Inequality~\eqref{ineq-main},
we arrive at
\[
x= \max_{1\leq k\leq K} \lambda_k + v_k \enspace .
\]
This shows that the family of vectors $\{v_k\}_{k=1,\ldots ,K}$ 
generates the tropical cone $\sC$.
Since the number of extreme rays of $\sC$ is bounded from above by the cardinality
of any of its generating families, this concludes the proof of Theorem~\ref{th-upperbound}.

\section{The tropical signed cyclic polyhedral cone and its polar}\label{SectionSigned}

We shall use the symmetrization of the max-plus semiring that M. Plus
introduced in~\cite{plus} to establish a max-plus analogue of the Cramer theorem. An intimately related Cramer theorem was established by Richter-Gebert, Sturmfels, and Theobald in~\cite{RGST}. In a nutshell, the result of~\cite{plus} deals
with max-plus linear systems in which signs are taken into account, whereas the
result of~\cite{RGST} concerns systems of equations in the tropical sense:
rather than requiring the maximum of ``positive'' terms of an expression to be equal to the maximum of its ``negative'' terms, it is only required that the maximum of the terms to be attained at least twice. The former Cramer theorem yields some
information on amoebas (image by the valuation) of linear spaces over the
field of real Puiseux series, whereas the latter Cramer theorem concerns amoebas over the field of complex Puiseux series. We refer the reader to the work by Akian, Gaubert and Guterman~\cite{AGG08b}, which gives a unified view of these Cramer theorems, connecting them also with a further work of Izhakian~\cite{Izh}.
In what follows, we need the version with signs, and use therefore
the result of~\cite{plus}, referring the reader to~\cite{AGG08b} for
more information.

The symmetrized max-plus semiring $\smax$ consists of three
copies of $\rmax$, glued by identifying the zero element.
A number of $\smax$ is written formally either as $a$,
$\ominus a$, or $a^\bullet$ for some $a\in \rmax$. These three
numbers are different, unless $a$ is the zero element (i.e.\ $a=-\infty$).
The {\em sign} $\sgn x$ of an element $x\in \smax $ is defined to be $+1$
if $x=a$ for some $a\in \rmax\setminus\{-\infty\}$, $-1$ if $x=\ominus a$
for some $a\in \rmax\setminus\{-\infty\}$, and $0$ otherwise. The elements
of the form $a$, $\ominus a$ and $a^\bullet$ are said to be
{\em positive}, {\em negative} and {\em balanced}, respectively. 
The elements which are either
positive or negative are said to be {\em signed}. We denote by $\smax^\vee$ the
set of signed elements and by $\smax^\bullet$ the set of balanced elements, so that 
\[
\smax =\smax^\vee \cup \smax^\bullet \enspace,
\]
the intersection of the latter sets being reduced to the zero element.
A vector is {\em signed} (resp.\ {\em balanced}) if each of its entries
is signed (resp.\ balanced).

The {\em modulus} of $x\in \{a,\ominus a,a^\bullet\}$
is defined as $|x|:=a$. 
The addition of two elements $x,y\in \smax$ is defined to
be $\max(|x|,|y|)$ if the maximum is attained only by elements
of positive sign, $\ominus \max(|x|,|y|)$ if it is attained
only by elements of negative sign, and $\max(|x|,|y|)^\bullet$
otherwise. 
For instance, $(\ominus 3)\oplus(2\oplus (\ominus 2))
= (\ominus 3) \oplus 2^\bullet=\ominus 3$. 
The multiplication
is defined in such a way that the modulus and the sign
are both morphisms. For instance, $(\ominus 3)(\ominus 4)=7$,
but $(\ominus 3)4^\bullet=7^\bullet$. 
The semiring $\smax$ is equipped with an involution
$x\mapsto \ominus x$, which sends the element $a$ to $\ominus a$, and vice
versa, and which fixes every balanced element $a^\bullet$. It is convenient to write,
for $x,y\in \smax$, $x\ominus y:=x\oplus (\ominus y)$.
We shall identify an element
$a\in \rmax$ with the corresponding element of $\smax$,
which yields an embedding of the semiring $\rmax$ into
$\smax$. 

The additive and multiplicative rules of $\smax$ become intuitive if the element
$a\in \smax$ is interpreted as the equivalence
class of real functions of $t$ belonging to $\Theta(t^a)$
as $t\to\infty$ (i.e., functions of $t$ belonging
to some interval $[Ct^a,C't^a]$ for some $C,C'>0$).
The element $\ominus a$ can be interpreted as the opposite of the
latter equivalence class, whereas $a^\bullet$ represents
the equivalence class $O(t^a)$. Then, the rule $(\ominus 3)\oplus(2\ominus 2)
= (\ominus 3) \oplus 2^\bullet=\ominus 3$ can be interpreted as the ``classical'' rule with asymptotic expansions: $-\Theta(t^3)+\Theta(t^2)-\Theta(t^2)=
-\Theta(t^3)+O(t^2)=-\Theta(t^3)$.

Given  $p$ scalars $-\infty<t_1< t_2< \cdots < t_p$ in $\rmax$,
and a collection of signs $\epsilon_{ij}\in \{\oplus \unit , \ominus \unit \}$, 
$1\leq i\leq p,1\leq j\leq d$, we construct the $p\times d$ matrix 
$C:=C(\epsilon,t)$ with entries in the symmetrized max-plus semiring
\[
C_{ij}=\epsilon_{ij} t_i^{j-1},\qquad 1\leq i\leq p\enspace ,\qquad 1\leq j\leq d 
\enspace .
\]
We denote by $C_i$ the $i$th row of $C$, 
and we write $C_i=C_i^+\ominus C_i^-$ where $C_i^+,C_i^-\in \rmax^{d}$ 
are chosen in such a way that for all $1\leq j\leq d$ 
exactly one of the $j$th entries of $C^+_i$ and $C^-_i$ is non-zero. 
\begin{definition}[Signed cyclic polyhedral cone]
The {\em signed cyclic polyhedral cone} with sign pattern 
$(\epsilon_{ij})$ is the tropical cone of $(\rmax^d)^2$ 
generated by the elements $(C^+_i,C^-_i)$, $1\leq i\leq p$. 
The {\em polar} of this cone is the set $\sK(\epsilon)$
of vectors 
$x\in \rmax^{d}$ such that 
\[
C_i^- x \leq C^+_{i}x \enspace  ,\qquad \forall 1\leq i\leq p \enspace .
\]
\end{definition}
The notion of tropical polar was introduced in~\cite{katz08}, to which
we refer the reader for more information.
We note that a related cyclic polytope (without signs) was studied
in by Block and Yu~\cite{blockyu06},

We shall often write $\sK$ instead of $\sK(\epsilon)$ for brevity.

We shall give a combinatorial construction of the extreme rays of $\sK$. 
An inequality $a x \leq b x$ ($a,b \in \rmax^d$) is said to be \emph{saturated} 
by $y \in \rmax^d$ if the equality $a y = b y$ holds. 
By analogy with the classical case, we expect an extreme generator to be
obtained by saturating $k$ inequalities among $C^+_{i}x\geq C_i^- x$,
$ 1\leq i\leq p$, and by setting $d-k-1$ entries of $x$ to zero.
In this way, we get $k$ equations for $k+1$ degrees of freedom, and
can hope the solution $x$ to be unique up to a scalar multiple.

In order to implement this idea, given two sequences of indices 
$I=\left\{i_1,\ldots ,i_k \right\}$ and 
$J=\left\{j_1,\ldots ,j_{k+1} \right\}$, where $k\leq d-1$ 
and $i_1<\cdots<i_k$, $j_1<\cdots<j_{k+1}$,
we consider the matrix $C(I,J)$ obtained by deleting the rows and columns 
of $C$ whose indices do not belong to $I$ and $J$, respectively. 
The matrices $C^+(I,J)$ and $C^-(I,J)$ are defined similarly.

We shall need to characterize the solutions $z$ of the system
$C^+(I,J)z=C^-(I,J)z$.

To this end, let us recall some basic consequences of the Cramer theorem
of~\cite{plus}. 
This result applies to systems of ``balances''. The balance
relation in $\smax$ is defined by $x\nabla y $ if $x\ominus y\in \smax^\bullet$. It is a non-transitive relation, which allows one to make
elimination arguments which are somehow similar to the case
of rings, although the addition does not have an opposite law.
In particular, if $x,y\in \smax$, $x=y$ implies that $x\ominus y\nabla\zero$,
and the converse holds if $x$ and $y$ are signed.
The balance relation is extended to vectors of $\smax^d$, being understood
entrywise.

Consider a linear system of the form
$A'x\oplus b'=A''x\oplus b''$, 
where $A',A''$ are $n\times n$ matrices with entries
in $\rmax$, and $b',b''\in \rmax^n$. Let
$A:=A'\ominus A''$, which is a well defined matrix
with entries in $\smax$. Similarly, let $b:=b''\ominus b'$.
It follows from the previous discussion that if
 $A'x\oplus b'=A''x\oplus b''$, then, the balance relation
$Ax\nabla b$ holds. Conversely, if $x$ is a vector with positive
entries, and if $Ax\nabla b$, then $A'x\oplus b'=A''x\oplus b''$.

The determinant of an $n\times n$ matrix $A=(a_{ij})$ 
with entries in $\smax $ is given by 
\[
\det A :=\bigoplus_{\sigma \in S_n} \sgn(\sigma ) 
a_{\sigma(1)1}\cdots a_{\sigma(n)n}
 \enspace ,  
\]
where $\sgn(\sigma ):=\oplus \unit $ if $\sigma $ is even and  
$\sgn(\sigma ):= \ominus \unit $ if $\sigma $ is odd. 
We denote by $A^{\adj }$ the transpose of the matrix of cofactors. 

\begin{theorem}[\cite{plus}, see also~{\cite[Th.~6.4]{AGG08b}}]\label{th-cramer}
Let $A$ be an $n\times n$ matrix with entries in $\smax $ and $b\in \smax^d$. 
Then, every signed solution of the system of balances 
\begin{equation}\label{SolutionBalance}
Ax\nabla b
\end{equation}
satisfies 
\[
\det A \; x \nabla A^{\adj } b \; .
\]
Conversely, if $A^{\adj } b$ is signed and if $\det A$ is invertible,
then $x=(\det A)^{-1} A^{\adj } b$ is the unique 
signed solution of~\eqref{SolutionBalance}. 
\end{theorem}  

By taking $b$ to be the zero vector, it follows that the equation $A x\nabla \zero$ has a non-zero signed solution only if $\det A$ is balanced. 
The converse implication also holds~\cite{plus}, but we shall not need
it here.

This max-plus analogue of Cramer theorem shows that the system 
of balances $Ax\nabla b$ can be solved by the usual Cramer rule, 
the determinants being interpreted as elements of $\smax$.
In particular, it shows that if none of the Cramer determinants
(the determinants appearing in the Cramer formula) is balanced,
then the system $Ax\nabla b$ has a unique signed solution, 
given by the Cramer formul\ae. Under the same circumstances, 
the original system $A'x\oplus b'=A''x\oplus b''$ has a solution in 
$\rmax$ if and only if the solution of the system of balances is positive.

We now apply this result to the homogeneous system $C^+(I,J)z=C^-(I,J)z$,
with $k$ equations and $k+1$ unknowns. 

Let us now consider the system of balances $C(I,J)z\nabla \zero$,
and let $D_r$ denote the $r$th Cramer determinant of this system,
which is the determinant of the matrix obtained from $C(I,J)$ by
deleting column $r$, i.e. $C(I,J\setminus \{ r\} )$.
\begin{lemma}\label{LemmaCramerDet}
The Cramer determinants of the previous linear system are given by
\begin{align*}
D_{k+1}&= t_{i_1}^{j_1-1}t_{i_2}^{j_2-1}\cdots t_{i_{k}}^{j_k-1}
\epsilon_{i_1 j_1}\epsilon_{i_2 j_2}\cdots \epsilon_{i_k j_k} 
\enspace , \\
D_1 &=t_{i_1}^{j_2-1}t_{i_2}^{j_3-1}\cdots t_{i_k}^{j_{k+1}-1}
\epsilon_{i_1 j_2}\epsilon_{i_2 j_3}\cdots \epsilon_{i_k j_{k+1}} \enspace , \\
D_r &= t_{i_1}^{j_1-1} \cdots t_{i_{r-1}}^{j_{r-1}-1} 
t_{i_r}^{j_{r+1}-1}\cdots t_{i_k}^{j_{k+1}-1} \epsilon_{i_1 j_1}\cdots 
\epsilon_{i_{r-1} j_{r-1}}\epsilon_{i_r j_{r+1}}\cdots 
\epsilon_{i_k j_{k+1}} \enspace , \enspace  2\leq r\leq k \enspace .
\end{align*}
\end{lemma}
\begin{proof}
When $A=C(I,J\setminus \{ r\} )$, we have 
\begin{equation}~\label{EqDet}
\det A =\bigoplus_{\sigma \in S_k} \sgn(\sigma ) 
\epsilon_{i_{\sigma(1)} j_1} t_{i_{\sigma(1)}}^{j_1-1} \cdots 
\epsilon_{i_{\sigma(r-1)} j_{r-1}} t_{i_{\sigma(r-1)}}^{j_{r-1}-1} 
\epsilon_{i_{\sigma(r)} j_{r+1}} t_{i_{\sigma(r)}}^{j_{r+1}-1}\cdots 
\epsilon_{i_{\sigma(k)} j_{k+1}} t_{i_{\sigma(k)}}^{j_{k+1}-1}
 \enspace , 
\end{equation}
since
\[
a_{\sigma (s) s}= 
\left\{
\begin{array}{ll}
\epsilon_{i_{\sigma (s)} j_s} t_{i_{\sigma (s)}}^{j_s-1} & \text{ if } s<r \enspace , \\
\epsilon_{i_{\sigma (s)} j_{s+1}} t_{i_{\sigma (s)}}^{j_{s+1}-1} & \text{ if } s\geq r \enspace . 
\end{array} 
\right. 
\]
If we define $\bar{\sigma }$ by $\bar{\sigma }(s)=s$ for $s=1,\ldots ,k$, 
it follows that 
\[ 
t_{i_{\sigma (1)}}^{j_1-1}\cdots  t_{i_{\sigma (r-1)}}^{j_{r-1}-1}t_{i_{\sigma (r)}}^{j_{r+1}-1}\cdots t_{i_{\sigma (k)}}^{j_{k+1}-1} < 
t_{i_{\bar{\sigma }(1)}}^{j_1-1}\cdots  t_{i_{\bar{\sigma }(r-1)}}^{j_{r-1}-1}
t_{i_{\bar{\sigma }(r)}}^{j_{r+1}-1}\cdots t_{i_{\bar{\sigma }(k)}}^{j_{k+1}-1} \enspace ,
\] 
for all $\sigma \neq \bar{\sigma }$, because 
$t_{i}^{j-1} t_{i'}^{j'-1}<t_{i'}^{j-1} t_{i}^{j'-1}$ whenever $t_{i'}<t_i$ and $j<j'$. 
Therefore, the term corresponding to $\bar{\sigma}$ in the only one 
maximizing the modulus in~\eqref{EqDet}, which implies that  
\[
D_r=\epsilon_{i_1 j_1} t_{i_1}^{j_1-1} \cdots \epsilon_{i_{r-1} j_{r-1}} t_{i_{r-1}}^{j_{r-1}-1} 
\epsilon_{i_r j_{r+1}} t_{i_r}^{j_{r+1}-1}\cdots \epsilon_{i_k j_{k+1}} t_{i_k}^{j_{k+1}-1} \enspace ,
\] 
and in particular $D_r$ is signed. 
\end{proof}
\begin{coro}\label{CoroSolBalance}
The system of balances $C(I,J)z\nabla \zero$ has a signed non-zero solution $z$,
which is unique up to a scalar multiple, and which is determined by the relations
\begin{equation}\label{RelationSol}
\begin{array}{lcl}
z_1 &=& \ominus t_{i_1}^{j_2-j_1} \epsilon_{i_1 j_1} \epsilon_{i_1 j_2} z_2 \\ 
z_2 &=& \ominus t_{i_2}^{j_3-j_2} \epsilon_{i_2 j_2} \epsilon_{i_2 j_3} z_3 \\ 
&\vdots&    \\ 
z_k &=& \ominus t_{i_k}^{j_{k+1}-j_k} \epsilon_{i_k j_k} 
\epsilon_{i_k j_{k+1}}z_{k+1} 
\end{array}
\end{equation}
\end{coro}
\begin{proof}
Let $A$ denote the matrix consisting of the first $k$ columns of
$C(I,J)$ and let $b$ denote the opposite of the last column of $C(I,J)$. 
Define $\bar{z}$ to be the vector consisting of the first $k$ coordinates of $z$.
Then, we have $C(I,J)z\nabla \zero$ if, and only if, $A\bar{z}\nabla bz_{k+1}$. 

The Cramer theorem above implies that 
\[ D_{k+1}z_r\nabla  (\ominus \unit)^{k-r+1} D_{r} z_{k+1},\qquad 1\leq r\leq k
\enspace .
\]
Recall that when two elements of $\smax$ $y$ and $y'$ are both signed,
$y\nabla y'$ implies $y=y'$. It follows that the relations~\eqref{RelationSol} hold. The same theorem also shows that, conversely, setting 
$z_{k+1}=\unit$, and defining $z$ by~\eqref{RelationSol}, we obtain
a solution of $C(I,J)z\nabla \zero$.
\end{proof}
We get as an immediate corollary.
\begin{coro}\label{CoroSol}
The system $C^+(I,J)z=C^-(I,J)z$ has a non-zero solution in $\rmax^n$ if and
only if 
\[\epsilon_{i_1 j_1} \epsilon_{i_1 j_2} = \epsilon_{i_2 j_2} \epsilon_{i_2 j_3} 
= \cdots =\epsilon_{i_k j_k} 
\epsilon_{i_k j_{k+1}} =\ominus \unit  \enspace.
\]
Then, this solution $z$ is determined by~\eqref{RelationSol}, up to a
scalar multiple. \qed
\end{coro}
We shall denote by $z(I,J)$ the vector defined by~\eqref{RelationSol} 
together with the normalization condition $z_{k+1}=\unit$. 
The vector $z(I,J)$ is a candidate to be an extreme generator of $\sK$. 
We shall see that only those subsets $I,J$ meeting a special combinatorial condition 
that we express in terms of lattice paths actually yield an extreme generator.

We shall visualize a pair of integers $(i,j)$, with
$1\leq i\leq p$ and $1\leq j\leq d$, as the position
of the corresponding entry in a $p\times d$ matrix.
So $(p,d)$ is the position of the bottom right entry
and $(1,1)$ is the position of the top left entry. 

\begin{figure}
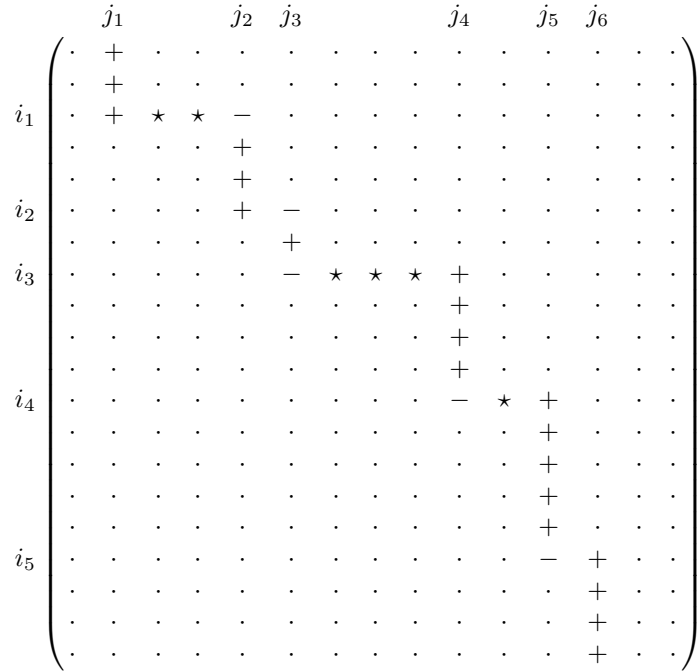

{\small
\[
\bordermatrix{
& & j_1 & & & j_2 & j_3 & & & & j_4 & & j_5 & j_6 & & \cr
&\cdot&\mop&\cdot&\cdot&\cdot&\cdot&\cdot&\cdot&\cdot&\cdot&\cdot&\cdot&\cdot&\cdot&\cdot \cr
&\cdot&\mop&\cdot&\cdot&\cdot&\cdot&\cdot&\cdot&\cdot&\cdot&\cdot&\cdot&\cdot&\cdot&\cdot \cr
i_1&\cdot& \mop&\star&\star&\mon&\cdot&\cdot&\cdot&\cdot&\cdot&\cdot&\cdot&\cdot&\cdot&\cdot\cr
&\cdot&\cdot&\cdot&\cdot&\mop&\cdot&\cdot&\cdot&\cdot&\cdot&\cdot&\cdot&\cdot&\cdot&\cdot\cr
&\cdot&\cdot&\cdot&\cdot&\mop&\cdot&\cdot&\cdot&\cdot&\cdot&\cdot&\cdot&\cdot&\cdot&\cdot\cr
i_2&\cdot&\cdot&\cdot&\cdot&\mop&\mon&\cdot&\cdot&\cdot&\cdot&\cdot&\cdot&\cdot&\cdot&\cdot\cr
&\cdot&\cdot&\cdot&\cdot&\cdot&\mop&\cdot&\cdot&\cdot&\cdot&\cdot&\cdot&\cdot&\cdot&\cdot\cr
i_3&\cdot&\cdot&\cdot&\cdot&\cdot&\mon&\star&\star&\star&\mop&\cdot&\cdot&\cdot&\cdot&\cdot\cr
&\cdot&\cdot&\cdot&\cdot&\cdot&\cdot&\cdot&\cdot&\cdot&\mop&\cdot&\cdot&\cdot&\cdot&\cdot\cr
&\cdot&\cdot&\cdot&\cdot&\cdot&\cdot&\cdot&\cdot&\cdot&\mop&\cdot&\cdot&\cdot&\cdot&\cdot\cr
&\cdot&\cdot&\cdot&\cdot&\cdot&\cdot&\cdot&\cdot&\cdot&\mop&\cdot&\cdot&\cdot&\cdot&\cdot\cr
i_4&\cdot&\cdot&\cdot&\cdot&\cdot&\cdot&\cdot&\cdot&\cdot&\mon&\star&\mop&\cdot&\cdot&\cdot\cr
&\cdot&\cdot&\cdot&\cdot&\cdot&\cdot&\cdot&\cdot&\cdot&\cdot&\cdot&\mop&\cdot&\cdot&\cdot\cr
&\cdot&\cdot&\cdot&\cdot&\cdot&\cdot&\cdot&\cdot&\cdot&\cdot&\cdot&\mop&\cdot&\cdot&\cdot\cr
&\cdot&\cdot&\cdot&\cdot&\cdot&\cdot&\cdot&\cdot&\cdot&\cdot&\cdot&\mop&\cdot&\cdot&\cdot\cr
&\cdot&\cdot&\cdot&\cdot&\cdot&\cdot&\cdot&\cdot&\cdot&\cdot&\cdot&\mop&\cdot&\cdot&\cdot\cr
i_5&\cdot&\cdot&\cdot&\cdot&\cdot&\cdot&\cdot&\cdot&\cdot&\cdot&\cdot&\mon&\mop&\cdot&\cdot\cr
&\cdot&\cdot&\cdot&\cdot&\cdot&\cdot&\cdot&\cdot&\cdot&\cdot&\cdot&\cdot&\mop&\cdot&\cdot\cr
&\cdot&\cdot&\cdot&\cdot&\cdot&\cdot&\cdot&\cdot&\cdot&\cdot&\cdot&\cdot&\mop&\cdot&\cdot\cr
&\cdot&\cdot&\cdot&\cdot&\cdot&\cdot&\cdot&\cdot&\cdot&\cdot&\cdot&\cdot&\mop&\cdot&\cdot\cr
}
\]}
\caption{A tropically allowed lattice path.}\label{fig:tropically_allowed_path}
\end{figure}

We shall consider oriented lattice paths, which are sequences of positions
starting from some top node $(1,j)$ and ending with some bottom node $(p,j)$,
in which at each step, the next position is either immediately at the right
or immediately at the bottom of the current one. Thus, 
such a path consists of vertical segments (oriented downward) 
and of horizontal segments (oriented from left to right). 
An example of lattice path is given in Figure~\ref{fig:tropically_allowed_path},
the initial (vertical) segment consists of the positions 
$(1,2)$, $(2,2)$, $(3,2)$, the next (horizontal) segment
consists of  $(3,2)$, $(3,3)$, $(3,4)$, $(3,5)$, 
the next (vertical) segment consists of 
$(3,5)$, $(4,5)$, $(5,5)$, $(6,5)$, etc. 
Note that the initial and final segments 
may be restricted to a unique position.

We shall say that such a lattice path is {\em tropically allowed} 
for the sign pattern $(\epsilon_{ij})$ if the following conditions are valid:
\begin{enumerate}[(i)]
\item\label{item:C1} every sign occurring on the initial vertical segment, 
except possibly the sign at the bottom of the segment, is positive;
\item\label{item:C2} every sign occurring on the final vertical segment, 
except possibly the sign at the top of the segment, is positive;
\item\label{item:C3} every sign occurring in some other vertical segment, 
except possibly the signs at the top and bottom of this segment, is positive;
\item\label{item:C4} for every horizontal segment, the pair of signs consisting 
of the signs of the leftmost and rightmost positions of the segment is of the
form $(+,-)$ or $(-,+)$; 
\item\label{item:C5} as soon as a pair $(-,+)$ occurs as the pair of 
extreme signs of some horizontal segment, 
the pairs of signs corresponding to all the horizontal segments 
below this one must also be equal to $(-,+)$.
\end{enumerate}
The notion of (non-tropically) allowed lattice path 
is defined only by Conditions~\eqref{item:C1}-\eqref{item:C4}. Hence, 
a tropically allowed path is allowed, 
but the converse is not true.

Figure~\ref{fig:tropically_allowed_path} gives an example of 
tropically allowed lattice path, 
the positions belonging to the path but the sign of which is irrelevant
are indicated by the symbol ``$\star$''. The positions which do not belong
to the path are indicated by the symbol ``$\cdot$''.

In order to prove Theorem~\ref{th-allowed},
it is convenient to recall the notion of tangent cone introduced in~\cite{AGG09}. 
Given a cone $\sC$ of $\rmax^d$ defined as the intersection of a finite set of 
tropical half-spaces $A_r x \leq B_r x$, 
where $A_r$ and $B_r$ denote the $r$th rows of some matrices $A$ and $B$,   
the \emph{tangent cone} of $\sC$ at $y\in \rmax^d$ is defined as the tropical 
cone $\sT(\sC,y)$ of $\rmax^d$ given by the system of inequalities
\begin{equation}
\max_{i \in \argmax(A_r y)} x_i \leq \max_{j \in \argmax(B_r y)} x_j \qquad 
\text{for all }r \text{ such that }A_r y = B_r y \enspace , 
\end{equation}
where $\argmax(c y)$ is the argument of the maximum 
$c y = \max_{1 \leq i \leq d} (c_i + y_i)$ for any row vector c. 
The tangent cone of $\sC$ at $y$ provides a local description of $\sC$ around $y$, 
leading to the following characterization of the extreme vectors of a tropical polyhedral cone. 

\begin{theorem}[\cite{AGG09}]\label{CharactExtreme} 
A vector $y\in \rmax^d$ belongs to an extreme ray of a tropical polyhedral cone $\sC$ if, 
and only if, there exists $s\in\{1,\dots ,d\}$ such that 
\begin{equation}\label{ExtremeProp}
(x\in \sT(\sC,y)\cap \{ \unit ,\zero \}^d \text{ and } x_s=\unit ) \implies
( x_r=\unit \text{ or } y_r=\zero )
\end{equation}
for all $r\in\{1,\ldots ,d\}$. 
\end{theorem}

As a consequence, we obtain. 
 
\begin{coro}\label{NumIneqSat}
Let $\sC: =\left\{x\in \rmax^d \mid A_r x\leq B_r x \; ,1\leq r \leq p \right\}$ 
be a tropical polyhedral cone and let $y\in \rmax^d$ 
be a vector in an extreme ray of $\sC$.  
If $t$ entries of $y$ are zero, 
then $y$ must saturate at least $d-t-1$ inequalities among 
$A_r x\leq B_r x$, $1\leq r \leq p$.  
\end{coro}

\begin{proof}
Let $s$ be an index satisfying the 
condition in Theorem~\ref{CharactExtreme}. 
Among the inequalities that define $\sT(\sC,y)$, 
consider those with precisely one term on the right hand side, 
i.e. those of the form 
\[ 
\oplus_{i\in I_h}x_i \leq x_h \enspace , 
\]
for some set of indices $I_h$. Let $H$ be the set composed of such indices $h$. 
If $y$ saturates strictly less than $d-t-1$ inequalities,  
there exists $q\in\{1,\ldots ,d\}$ such that  
$q \not \in \{s\}\cup H \cup \left\{ j\mid y_j=\zero \right\}$. 
Then, the vector $x\ \in \{ \unit ,\zero \}^d$ defined by $x_q:=\zero$ 
and $x_i:=\unit $ for all $i\neq q$ belongs to $\sT(\sK,y)$, 
which contradicts~\eqref{ExtremeProp}. 
\end{proof}

We now prove Theorem~\ref{th-allowed}.
Let $x\in \rmax^d$ be a vector in an extreme ray of $\sK$. 
Assume that $x_j\neq \zero$ if, and only if, 
$j\in J=\left\{j_1,\ldots ,j_{k+1} \right\}$, where $k\leq d-1$.
Then, by Corollary~\ref{NumIneqSat} we know 
that $x$ must saturate at least $k$ 
inequalities among $C_i^- x\leq C_i^+ x$, $i=1,\ldots ,p$. 
More precisely, we claim that $x$ saturates exactly $k$ inequalities. 
To see this, let $\bar{x}\in \rmax^{k+1}$ be the vector obtained by 
deleting the entries of $x$ which do not belong to $J$ 
(or equivalently, are zero).
Assume that $x$ saturates the inequalities $C_i^- x\leq C_i^+ x$ for $i\in I$, 
where $I$ has $k+1$ elements. Then, we would have $C(I,J)\bar{x}\nabla \zero$, 
where the determinant of the $(k+1)\times (k+1)$ matrix 
$C(I,J)$ is signed by Lemma~\ref{LemmaCramerDet}, 
contradicting Cramer theorem above in the 
case of homogeneous systems of balances. This proves our claim.  

Let $I=\left\{i_1,\ldots ,i_k \right\}$ be the set composed of the indices of 
the inequalities which $x$ saturates. We assume  
$i_1<\cdots <i_k$ and $j_1<\cdots <j_{k+1}$.  
With this extreme ray, we associate the lattice path $\sP$
\begin{equation}\label{Path}
(1,j_1),\ldots ,(i_1,j_1),\ldots ,(i_1,j_2),\ldots , (i_2,j_2),
\ldots ,(i_k ,j_k),\ldots ,(i_k,j_{k+1}),\dots ,(p,j_{k+1}) \; .
\end{equation} 
In other words, the ordinates of the horizontal segments of this path are 
given by the indices of the inequalities which are saturated, 
and the absciss\ae\ of the vertical segments are given by the indices 
$j$ such that $x_j$ is non-zero. 
Note that this path has $k$ horizontal segments and that 
$(i_r,j_r)$ and $(i_r,j_{r+1})$ are the leftmost and rightmost 
positions of the $r$th horizontal segment.  

We claim that $\sP$ is tropically allowed. In order to prove this, 
define $\bar{x}\in \rmax^{k+1}$ as above.  
Since $C^+(I,J)\bar{x} =C^-(I,J)\bar{x} $, 
the ``only if'' part of Corollary~\ref{CoroSol}  
shows precisely that $\sP$ satisfies Condition~\eqref{item:C4}.  
Hence, we may assume  
\begin{equation}\label{sol}
\bar{x} = z(I,J) = 
\begin{pmatrix}
t_{i_1}^{j_2-j_1}  t_{i_2}^{j_3-j_2}  \ldots t_{i_k}^{j_{k+1}-j_k}  \cr
t_{i_2}^{j_3-j_2} \ldots t_{i_k}^{j_{k+1}-j_k} \cr
\vdots \cr
t_{i_k}^{j_{k+1}-j_k} \cr
\unit
\end{pmatrix} \; , 
\end{equation}
which implies that
\begin{equation}\label{Cix}
(C_i^+\oplus C_i^-) x = 
\bigoplus_{1\leq j \leq d} t_i^{j-1} x_j = 
\bigoplus_{1\leq r\leq k+1} t_i^{j_r-1} x_{j_r} =
\bigoplus_{1\leq r\leq k+1} t_i^{j_r-1} t_{i_r}^{j_{r+1}-j_r} 
\ldots t_{i_k}^{j_{k+1}-j_k} \enspace ,
\end{equation}
for all $1\leq i \leq p$. Then, if $i_s < i < i_{s+1}$, it follows that 
\[
t_i^{j_r-1} t_{i_r}^{j_{r+1}-j_r}\ldots t_{i_k}^{j_{k+1}-j_k} < 
t_i^{j_{s+1}-1} t_{i_{s+1}}^{j_{s+2}-j_{s+1}}\ldots t_{i_k}^{j_{k+1}-j_k} 
\enspace ,
\] 
for all $r\neq s+1$. Hence, the maximum in~\eqref{Cix} is attained 
only for $j=j_{s+1}$. Since $C_i^- x\leq C_i^+ x$, 
the sign of $C_{i j_{s+1}}$ must be positive, 
implying that Condition~\eqref{item:C3} is valid for $\sP$. 
Analogously, when $i<i_1$
\[
t_i^{j_r-1} t_{i_r}^{j_{r+1}-j_r} \ldots t_{i_k}^{j_{k+1}-j_k} 
< t_i^{j_1-1} t_{i_1}^{j_2-j_1} \ldots t_{i_k}^{j_{k+1}-j_k} 
\enspace ,
\] 
for all $r>1$, so $C_i^- x\leq C_i^+ x$ implies that Condition~\eqref{item:C1} 
holds for $\sP$. Finally, if $i>i_k$, we have 
\[
t_i^{j_r-1} t_{i_r}^{j_{r+1}-j_r} \ldots t_{i_k}^{j_{k+1}-j_k} 
< t_i^{j_{k+1}-1}  
\enspace ,
\] 
for all $r<k+1$, and the same argument as before shows that  
$\sP$ satisfies Condition~\eqref{item:C2}. 

When $i=i_s$ for some $s$, we have
\[
t_i^{j_r-1} t_{i_r}^{j_{r+1}-j_r} \ldots t_{i_k}^{j_{k+1}-j_k} < 
t_i^{j_s-1} t_{i_s}^{j_{s+1}-j_s}\ldots t_{i_k}^{j_{k+1}-j_k} =
t_i^{j_{s+1}-1} t_{i_{s+1}}^{j_{s+2}-j_{s+1}}\ldots t_{i_k}^{j_{k+1}-j_k} 
\enspace ,
\]
for all $r\not \in \left\{s ,s+1\right\}$, 
which means that the tangent cone $\sT(\sK,x)$ of $\sK$ at 
$x$ is defined by the inequalities 
$x_{j_r}\geq x_{j_{r+1}}$ if 
$(\epsilon_{i_r j_r},\epsilon_{i_r j_{r+1}})=(\oplus \unit, \ominus \unit)$ 
and $x_{j_r}\leq x_{j_{r+1}}$ if
$(\epsilon_{i_r j_r},\epsilon_{i_r j_{r+1}})=(\ominus \unit, \oplus \unit)$, 
for $r=1,\ldots ,k$. 
It is convenient to visualize the relations defining the tangent 
cone by constructing a digraph with nodes $j_{1}, \ldots ,j_{k+1}$ 
and an arc from $j_r$ to $j_{r+1}$ (resp.\ from $j_{r+1}$ to $j_r$) 
if the inequality $x_{j_r}\geq x_{j_{r+1}}$ (resp.\ $x_{j_{r+1}}\geq x_{j_r}$)  
belongs to these relations. 
For instance, the digraph associated with the relations 
$x_{j_1}\geq x_{j_2}\geq x_{j_3}\leq x_{j_4}\leq x_{j_5} \geq x_{j_6}$ is 
\[
j_1 \rightarrow j_2 \rightarrow j_3 \leftarrow 
j_4 \leftarrow j_5 \rightarrow j_6 
\] 
Theorem~\ref{CharactExtreme} requires the existence of a node 
$j_s$ such that $x_{j_s}=\unit$ implies $x_{j_r}=\unit$ for all $r$. 
This can only occur if in the digraph associated with the tangent cone 
there is a directed path from any node to $j_s$.
Since the digraph associated with $\sT(\sK,x)$ has a line structure,
the only possibility for this to happen is that, 
when scanning the arcs of the digraph from left to right, 
the arcs must be directed to the right until node $j_s$, 
and then all the remaining arcs must be directed to the left. 
Since an arc directed to the right (resp.\ left) 
corresponds to an horizontal segment of the path 
whose pair of extreme signs is $(+,-)$ (resp.\ $(-,+)$),  
it follows that $\sP$ must satisfy Condition~\eqref{item:C5}. 
In consequence, $\sP$ is tropically allowed. 

Conversely, with a tropically allowed lattice path with $k$ 
horizontal segments, we associate the sequences of indices 
$I=\left\{i_1,\ldots,i_k\right\}$ and $J=\left\{ j_1,\dots ,j_{k+1}\right\}$ 
obtained by taking the ordinates and absciss\ae\ of its 
horizontal and vertical segments respectively, 
as illustrated for the tropically allowed path 
in Figure~\ref{fig:tropically_allowed_path}. 
Note that the previous logic is reversible, meaning that if we 
define $x\in \rmax^d$ by $x_{j_r}=z_r(I,J)$ for 
$r\in \left\{1,\ldots ,k+1\right\}$ and 
$x_j=\zero $ for $j\not \in J$, then $x$ is in an extreme ray of $\sK$. 
More precisely, by Corollary~\ref{CoroSol}, Condition~\eqref{item:C4} implies 
that all the entries of $z(I,J)$ are positive so $x\in \rmax^d$, 
Conditions~\eqref{item:C1},~\eqref{item:C2} and~\eqref{item:C3} imply that $x$ belongs to $\sK$, 
and finally Condition~\eqref{item:C5} and Theorem~\ref{CharactExtreme} 
show that $x$ belongs to an extreme ray of $\sK$.   
This concludes the proof of Theorem~\ref{th-allowed}.

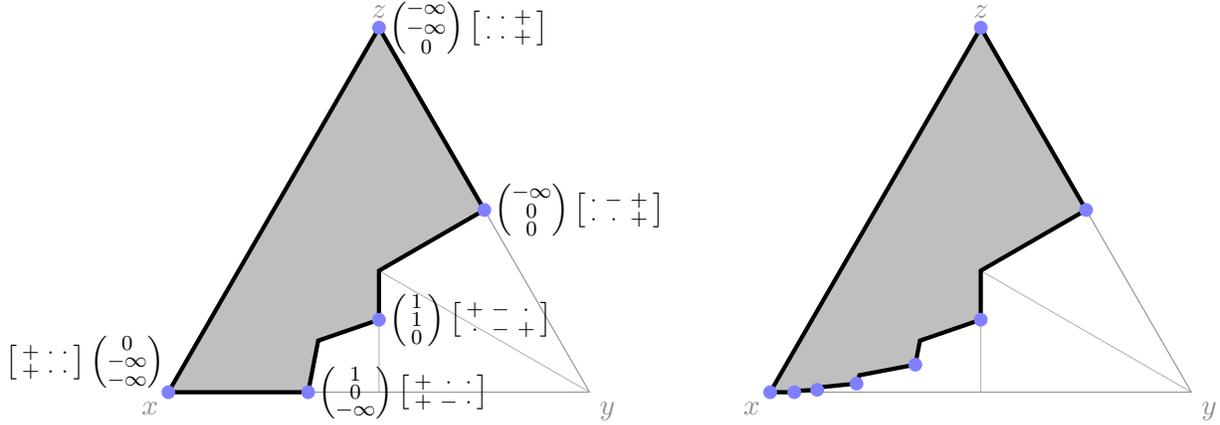
\begin{figure}[t]
\begin{center}
\begin{tikzpicture}[scale=0.8,convex/.style={draw=lightgray,fill=lightgray},point/.style={blue!50},line/.style={blue!50,ultra thick},convexborder/.style={ultra thick},pointlabel/.style={black}]
\equilateral{7}{90};

\barycenter{g0}{\expo{0}}{\expo{-1}}{0};
\barycenter{g1}{\expo{1}}{\expo{1}}{\expo{0}};
\barycenter{g2}{0}{\expo{0}}{\expo{0}};

\barycenter{g01}{\expo{2}}{\expo{1}}{\expo{0}};
\barycenter{g12}{\expo{0}}{\expo{0}}{\expo{0}};

\filldraw[convex] (x) -- (g0) -- (g01) -- (g1) -- (g12) -- (g2) -- (z) -- cycle;
\draw[convexborder] (x) -- (g0) -- (g01) -- (g1) -- (g12) -- (g2) -- (z) -- cycle;

\filldraw[point] (x) circle (3pt) node[black,above left=-0.5ex] 
{$\left[\begin{smallmatrix} + & \cdot & \cdot \\ + & \cdot & \cdot \end{smallmatrix}\right]
\left(\!\begin{smallmatrix}0 \\ -\infty \\ -\infty\end{smallmatrix}\!\right)$};
\filldraw[point] (z) circle (3pt) node[black,right] 
{$\left(\!\begin{smallmatrix}-\infty \\ -\infty \\ 0\end{smallmatrix}\!\right)%
\left[\begin{smallmatrix} \cdot & \cdot & + \\ \cdot & \cdot & +  \end{smallmatrix}\right]$};
\filldraw[point] (g0) circle (3pt) node[black,right] 
{$\left(\!\begin{smallmatrix}1 \\ 0 \\ -\infty\end{smallmatrix}\!\right)%
\left[\begin{smallmatrix} + & \cdot & \cdot  \\ + & - & \cdot \end{smallmatrix}\right]$};
\filldraw[point] (g1) circle (3pt) node[black,right] 
{$
\left(\begin{smallmatrix}1 \\ 1 \\ 0\end{smallmatrix}\right)%
\left[\begin{smallmatrix} + & - & \cdot \\ \cdot & - & +  \end{smallmatrix}\right]$};
\filldraw[point] (g2) circle (3pt) node[black,right] 
{$\left(\!\begin{smallmatrix}-\infty \\ 0 \\ 0\end{smallmatrix}\!\right) 
\left[\begin{smallmatrix} \cdot & - & + \\ \cdot & \cdot & +  \end{smallmatrix}\right]$};

\begin{scope}[xshift=10cm]
\equilateral{7}{90};

\barycenter{g0}{\expo{0}}{\expo{-4}}{0};
\barycenter{g1}{\expo{7}}{\expo{4}}{\expo{0}};
\barycenter{g2}{\expo{5}}{\expo{3}}{\expo{0}};
\barycenter{g3}{\expo{3}}{\expo{2}}{\expo{0}};
\barycenter{g4}{\expo{1}}{\expo{1}}{\expo{0}};
\barycenter{g5}{0}{\expo{0}}{\expo{0}};

\barycenter{g01}{\expo{8}}{\expo{4}}{\expo{0}};
\barycenter{g12}{\expo{6}}{\expo{3}}{\expo{0}};
\barycenter{g23}{\expo{4}}{\expo{2}}{\expo{0}};
\barycenter{g34}{\expo{2}}{\expo{1}}{\expo{0}};
\barycenter{g45}{\expo{0}}{\expo{0}}{\expo{0}};

\filldraw[convex] (x) -- (g0) -- (g01) -- (g1) -- (g12) -- (g2) -- (g23) -- (g3) -- (g34) -- (g4) -- (g45) -- (g5) -- (z) -- cycle;
\draw[convexborder] (x) -- (g0) -- (g01) -- (g1) -- (g12) -- (g2) -- (g23) -- (g3) -- (g34) -- (g4) -- (g45) -- (g5) -- (z) -- cycle;

\filldraw[point] (x) circle (3pt);
\filldraw[point] (z) circle (3pt);
\filldraw[point] (g0) circle (3pt);
\filldraw[point] (g1) circle (3pt);
\filldraw[point] (g2) circle (3pt);
\filldraw[point] (g3) circle (3pt);
\filldraw[point] (g4) circle (3pt);
\filldraw[point] (g5) circle (3pt);
\end{scope}
\end{tikzpicture}
\end{center}
\caption{The polars of two signed cyclic polyhedral cones in $\rmax^3$.}\label{fig:cones}
\end{figure}

\begin{example}
Figure~\ref{fig:cones} provides two examples of polars of signed cyclic polyhedral cones for $d = 3$. They are represented in barycentric coordinates: each element $(x_1,x_2,x_3)$ of $\rmax^3$ is represented as a barycenter with weights $(e^{x_1}, e^{x_2}, e^{x_3})$ of the three vertices of the outermost triangle. Then two representatives of a same ray are represented by the same point. This barycentric representation is convenient to represent points with infinite coordinates, which are mapped to the boundary of the simplex.

The two cones are defined by $p = 2$ and $p = 5$ inequalities respectively, and, for all $1 \leq i \leq p$, $t_i = i-1$ and $\epsilon_{ij} = \ominus \unit$ if and only if $j=2$. In other words, the first cone is associated with the sign pattern $\left(\begin{smallmatrix} + & - & + \\ + & - & + \end{smallmatrix}\right)$, and its polar is defined as the set of elements $(x_1,x_2,x_3) \in \rmax^3$ such that:
\[
\begin{pmatrix}
0 & -\infty & 0 \\
0 & -\infty & 2 
\end{pmatrix}
\begin{pmatrix}
x_1 \\ x_2 \\ x_3
\end{pmatrix}
\geq 
\begin{pmatrix}
-\infty & 0 & -\infty \\
-\infty & 1 & -\infty 
\end{pmatrix}
\begin{pmatrix}
x_1 \\ x_2 \\ x_3 
\end{pmatrix}.
\]

The extreme rays are depicted by blue points. For the first cone, a representative of each extreme ray is provided, and the corresponding tropically allowed path is given beside.
\end{example}

\section{The number of extreme points of the classical polar of a signed cyclic polyhedral cone}

We next give a characterization of the extreme rays of the polar of 
the classical analogue of the signed cyclic polyhedral cone, 
which shows that in the tropical case there exist fewer extreme rays. 
Therefore, in this section, all the operations should be 
understood in the usual algebra. 

Given $p$ positive real numbers $t_1<\cdots <t_p$   
and a sign pattern $(\epsilon_{ij})$, which now belongs 
to $\{+1,-1\}^{p\times d}$, 
we shall consider the usual polar of the signed cyclic polyhedral cone, 
which we still denote by $\sK$,
\[
\sK:=\left\{x\in\R^d \mid x\geq 0\; , C x \geq 0 \right\} \enspace ,
\]
where $C_{ij} =\epsilon_{ij} t_i^{j-1}$ 
for $1\leq i \leq p$ and $1\leq j \leq d$. 
 
Like in the previous section, given two sequences of indices 
$I=\{i_1,\ldots ,i_k \}$ and 
$J=\{j_1,\ldots ,j_{k+1}\}$ where $k\leq d-1$, 
consider the matrix $C(I,J)$ 
obtained by keeping only the rows and columns 
of $C$ whose indices belong to $I$ and $J$, respectively. 

\begin{lemma}\label{LemmaDetClass} 
The (ordinary) Cramer determinants $D_r$ 
of the system $C(I,J)z=0$ are given by
\begin{align}\label{FormDetCramer}
D_r =  \delta_r t_{i_1}^{j_1-1} \cdots t_{i_{r-1}}^{j_{r-1}-1} 
t_{i_r}^{j_{r+1}-1}\cdots t_{i_k}^{j_{k+1}-1} \epsilon_{i_1 j_1}\cdots 
\epsilon_{i_{r-1} j_{r-1}}\epsilon_{i_r j_{r+1}}\cdots 
\epsilon_{i_k j_{k+1}} \enspace , \enspace  
\end{align}
for $1\leq r\leq k+1$, where the scalars $\delta_r$ tend to $1$ as 
the ratios $t_2/t_1 , \ldots , t_{p+1}/t_p$ tend to infinity. 
\end{lemma} 

\begin{proof}
For any permutation $\sigma \in S_k$ define 
\[
D_r(\sigma ):=  
t_{i_{\sigma(1)}}^{j_1-1} \cdots t_{i_{\sigma(r-1)}}^{j_{r-1}-1} 
t_{i_{\sigma(r)}}^{j_{r+1}-1}\cdots t_{i_{\sigma(k)}}^{j_{k+1}-1}
 \enspace , 
\] 
so that 
\[
D_r =\sum_{\sigma \in S_k} \sgn(\sigma ) D_r(\sigma ) 
\epsilon_{i_{\sigma(1)} j_1} \cdots \epsilon_{i_{\sigma(r-1)} j_{r-1}} 
\epsilon_{i_{\sigma(r)} j_{r+1}} \cdots \epsilon_{i_{\sigma(k)} j_{k+1}} 
\enspace . 
\]
Let $\bar{\sigma}$ be defined by $\bar{\sigma}(s):=s$ for all $s$. 
We claim that for any $\sigma \neq \bar{\sigma}$, the quotient 
$D_r(\sigma )/D_r(\bar{\sigma })$ is a product of terms of the form 
$t_{i_r}/t_{i_s}$, where $i_r < i_s$. To see this, 
let $s=\max\{h\mid {\sigma(h)}\neq h\}$. Then, 
we have $\sigma (s) < s$ and there exists $q<s$ 
such that $s=\sigma (q)$. If we define $\sigma'$ by $\sigma'(s)=s$, 
$\sigma'(q)=\sigma(s)$ and $\sigma'(h)=\sigma (h)$ 
for all $h\not \in \{ s,q \}$, it follows that 
\[
\frac{D_r(\sigma )}{D_r(\sigma')} =
\frac{t_{i_{\sigma (s)}}^{j_{\hat{s}}-1} t_{i_s}^{j_{\hat{q}}-1}}
{t_{i_s}^{j_{\hat{s}}-1} t_{i_{\sigma (s)}}^{j_{\hat{q}}-1}} =
\left( \frac{t_{i_{\sigma (s)}}}{t_{i_s}} \right)^{j_{\hat{s}} - j_{\hat{q}}} 
\enspace , 
\]
where $\hat{s}=s$ if $s<r$ and $\hat{s}=s+1$ otherwise, 
and the same applies to $q$. 
The claim follows by repeating this procedure till $\sigma'=\bar{\sigma}$.  

Note that~\eqref{FormDetCramer} is satisfied for  
\begin{align*}
\delta_r := 1+\sum_{\sigma\neq \bar{\sigma }} 
\frac{\sgn(\sigma ) D_r(\sigma )\epsilon_{i_{\sigma(1)} j_1} \cdots 
\epsilon_{i_{\sigma(r-1)} j_{r-1}} \epsilon_{i_{\sigma(r)} j_{r+1}} 
\cdots \epsilon_{i_{\sigma(k)} j_{k+1}}}
{ D_r(\bar{\sigma }) \epsilon_{i_1 j_1}\cdots 
\epsilon_{i_{r-1} j_{r-1}}\epsilon_{i_r j_{r+1}}\cdots 
\epsilon_{i_k j_{k+1}}} \enspace ,   
\end{align*}
and from the discussion above it follows that $\delta_r$ tends to 
$1$ as the ratios $t_2/t_1 , \ldots , t_{p+1}/t_p$ tend to infinity.
\end{proof}

As a consequence of the classical Cramer theorem we obtain.

\begin{coro}\label{CoroSolClass}
Assume that the ratios $t_2/t_1 , \ldots , t_{p+1}/t_p$ 
are sufficiently large. Then, the system $C(I,J)z=0$ 
has a non-zero solution, which is unique up to a scalar multiple, 
and which is determined by the relations
\begin{equation}\label{RelationSolClassical}
\begin{array}{lcl}
z_1 &=& (-\gamma_1) t_{i_1}^{j_2-j_1} \epsilon_{i_1 j_1} 
\epsilon_{i_1 j_2} z_2 \enspace ,\\ 
z_2 &=& (-\gamma_2) t_{i_2}^{j_3-j_2} \epsilon_{i_2 j_2} 
\epsilon_{i_2 j_3} z_3 \enspace ,\\ 
&\vdots&    \\ 
z_k &=& (-\gamma_k) t_{i_k}^{j_{k+1}-j_k} \epsilon_{i_k j_k} 
\epsilon_{i_k j_{k+1}}z_{k+1} \enspace ,
\end{array}
\end{equation}
where for $r=1,\ldots ,k$ the scalars $\gamma_r$ tend to $1$ as the ratios 
$t_2/t_1 ,\ldots , t_{p+1}/t_p$ tend to infinity.  
\end{coro}

Like in the previous section, we shall denote by $z(I,J)$ the vector defined 
by~\eqref{RelationSolClassical} together with the normalization 
condition $z_{k+1}=1$, i.e. 
\begin{equation}\label{SolClass}
z(I,J): = 
\begin{pmatrix}
(-\gamma_1) \epsilon_{i_1 j_1} \epsilon_{i_1 j_2} 
  t_{i_1}^{j_2-j_1}  t_{i_2}^{j_3-j_2} \ldots t_{i_k}^{j_{k+1}-j_k} \cr
(-\gamma_2) \epsilon_{i_2 j_2} \epsilon_{i_2 j_3} 
t_{i_2}^{j_3-j_2} \ldots t_{i_k}^{j_{k+1}-j_k} \cr
\vdots \cr
(-\gamma_k) \epsilon_{i_k j_k} \epsilon_{i_k j_{k+1}} 
t_{i_k}^{j_{k+1}-j_k}  \cr
1
\end{pmatrix} \; . 
\end{equation} 

\begin{theorem} \label{th-class}
If the ratios $t_2/t_1 , \ldots , t_{p+1}/t_p$ are sufficiently large, 
the extreme rays of $\sK$ are in one to one correspondence 
with the (non-tropically) allowed lattice paths for the sign 
pattern $(\epsilon_{ij})$.   
\end{theorem}

\begin{proof}
Let $x\in \R^d$ be in an extreme ray of $\sK$. Assume that 
$\left\{j\mid x_j\neq 0\right\}=\left\{j_1,\dots ,j_{k+1}\right\}$, 
where $k\leq d-1$. Then, $x$ must saturate at least $k$ inequalities 
among $C_i x\geq 0$, $i=1,\ldots ,p$. Indeed, like in the tropical case, 
$x$ saturates precisely $k$ inequalities, because otherwise, 
by Lemma~\ref{LemmaDetClass} and Cramer theorem, 
it would be equal to the null vector. 

Let $\left\{i_1,\dots ,i_k\right\}$ be the indices of the inequalities 
which $x$ saturates. We assume $i_1<\cdots <i_k$ and $j_1<\dots <j_{k+1}$.  
With $x$ we associate the lattice path $\sP$ defined by~\eqref{Path}. 
We next show that $\sP$ is allowed. 

Let $\bar{x}\in \R^{k+1}$ be the vector obtained 
from $x$ by deleting its null entries. 
Since $\bar{x}$ satisfies $C(I,J)\bar{x}=0$
and the entries of $\bar{x}$ are positive, 
by Corollary~\ref{CoroSolClass} it follows that the signs on the 
extreme positions of every horizontal segment of $\sP$, 
i.e. $\epsilon_{i_r j_r}$ and $\epsilon_{i_r j_{r+1}}$, 
must be opposite. In other words, $\sP$ satisfies Condition~\eqref{item:C4}. 

Since we may assume $\bar{x}=z(I,J)$, it follows that 
\begin{equation}
C_i x = 
\sum_{1\leq j \leq d} \epsilon_{ij} t_i^{j-1} x_j = 
\sum_{1\leq r\leq k+1} \epsilon_{ij_r} t_i^{j_r-1} x_{j_r} =
\sum_{1\leq r\leq k+1} \epsilon_{ij_r} t_i^{j_r-1} 
\gamma_r t_{i_r}^{j_{r+1}-j_r} \ldots t_{i_k}^{j_{k+1}-j_k} \enspace ,
\end{equation}
for all $1\leq i \leq p$, where we define $\gamma_{k+1}:=1$.

If we take $i_s<i<i_{s+1}$, note that  
\[
C_i x=\epsilon_{ij_{s+1}} t_i^{j_{s+1}-1} \gamma_{s+1} 
t_{i_{s+1}}^{j_{s+2}-j_{s+1}} \ldots t_{i_k}^{j_{k+1}-j_k} 
(1+ \kappa_i) 
\]
where
\[
\kappa_i= \sum_{1\leq r \leq s} 
\frac{\gamma_r}{\gamma_{s+1}}
\left( \frac{t_{i_r}}{t_i} \right)^{j_{r+1} - j_r} \cdots 
\left( \frac{t_{i_s}}{t_i} \right)^{j_{s+1} - j_s} +
\sum_{s+2\leq r \leq k+1} 
\frac{\gamma_r}{\gamma_{s+1}}
\left( \frac{t_{i}}{t_{i_{s+2}}} \right)^{j_{s+3}-j_{s+2}} \cdots
\left( \frac{t_{i}}{t_{i_{r}}} \right)^{j_{r+1}-j_r} \enspace .
\]
Since $\kappa_i$ tend to $0$ as the ratios 
$t_2/t_1 , \ldots , t_{p+1}/t_p$ tend to infinity, 
it follows that $\epsilon_{ij_{s+1}}=+1$ must be satisfied 
in order to have $C_ix\geq 0$. 
This means that Condition~\eqref{item:C3} is valid for $\sP$. 
A similar argument shows that Conditions~\eqref{item:C1} and~\eqref{item:C2} also hold 
and thus $\sP$ is allowed. 

Conversely, like in the tropical case, 
note that the previous logic is reversible. 
With an allowed lattice path with $k$ horizontal segments, 
we associate the sequences of indices $I=\left\{i_1,\ldots,i_k\right\}$ 
and $J=\left\{ j_1,\dots ,j_{k+1}\right\}$ 
obtained by taking the ordinates and absciss\ae\ of 
its horizontal and vertical segments respectively.
If we define $x\in \R^d$ by $x_{j_r}=z_r(I,J)$ 
for $r\in \left\{1,\ldots ,k+1\right\}$ and 
$x_j=0$ for $j\not \in J$, then $x$ is in an extreme ray of 
$\sK$ for sufficiently large ratios $t_2/t_1 ,\ldots ,t_{p+1}/t_p$. 
More precisely, by Corollary~\ref{CoroSolClass}, Condition~\eqref{item:C4} implies 
that all the entries of $z(I,J)$ are positive, so $x_j\geq 0$ for all $j$. 
This fact together with Conditions~\eqref{item:C1},~\eqref{item:C2} and~\eqref{item:C3} 
imply that $x$ belongs to $\sK$.  
Finally, note that Lemma~\ref{LemmaDetClass}  
shows that the gradients of the inequalities that $x$ saturates, 
i.e.\ $C_{i_r}x\geq 0$ for $r\in I$ and $x_j\geq 0$ for $j\not \in J$, 
form a family of full rank. Therefore, 
$x$ belongs to an extreme ray of $\sK$. This concludes the proof. 
\end{proof}

The following theorem shows that
the bound $U(p+d,d-1)$ is attained by the polar of the classical 
analogue of the signed cyclic polyhedral cone. Its proof also shows
that the lattice path characterization of Theorem~\ref{th-class} may be 
thought of as a generalization of Gale's evenness criterion, since
the latter is recovered by considering the special case in which
the sign pattern is $\epsilon_{ij}=(-1)^j$.

\begin{theorem}\label{theo-classical}
The number of extreme rays of the classical polar $\sK$ of the 
signed cyclic polyhedral cone with sign pattern 
$\epsilon_{ij}:=(-1)^{j-1}$ is exactly $U(p+d,d-1)$. 
\end{theorem}

\begin{proof}
Given the set $\{1,\ldots ,n\}$, 
we shall say that a subset $Q$ of $\{1,\ldots ,n\}$ 
satisfies Gale's evenness condition, 
if for any $i,j\in \{1,\ldots ,n\}\setminus Q$ 
the number of elements in $Q$ between $i$ and $j$ is even. 
It is known (see~\cite{matousek}) 
that the number of subsets $Q$ of $\{1,\ldots ,n\}$ with $k$ elements 
satisfying the evenness condition is $U(n,k)$. 
We shall show that the number of extreme rays of $\sK$ is 
$U(p+d,d-1)$ by constructing a bijective correspondence between 
allowed lattice paths for the sign pattern $(\epsilon_{ij})$ and 
subsets of $\{1,\ldots ,p+q\}$ with $d-1$ elements which satisfy 
Gale's evenness condition. 

Given an allowed lattice path $\sP$ 
for the sign pattern $(\epsilon_{ij})$, 
let $I=\left\{i_1,\ldots,i_k\right\}$ and 
$J=\left\{ j_1,\dots ,j_{k+1}\right\}$ 
be the sets of ordinates and absciss\ae\ of 
its horizontal and vertical segments respectively. 
With $\sP$ we associate the subset $Q$ of 
$\{1,\ldots ,p+q\}$ defined by 
\[
Q :=\{i+d\mid i\in I\} \cup \{d-j+1\mid j\not \in J\} \enspace .
\] 
The set $Q$ may be visualized by scanning first the columns of the matrix
$\epsilon$ from right to left, keeping only the columns not in $J$,
and scanning then the rows of $\epsilon$ from top to bottom, keeping
now the rows in $I$. The following illustrates the definition of $Q$
for a special lattice path, the elements of $Q$ are listed by numbers on
the top and left borders of the matrix so that $Q=\{2,3,7,8,10,11,13,14\}$:
\[
\def\bul{\bullet}
\bordermatrix{&\cdot &8 &7 &\cdot & \cdot&\cdot & 3&2&\cdot\cr
10 &+&-&+&-&\cdot&\cdot&\cdot&\cdot&\cdot\cr
11 &\cdot&\cdot &\cdot&- &+&\cdot&\cdot&\cdot&\cdot\cr
\,\,\cdot&\cdot &\cdot &\cdot &\cdot &+&\cdot&\cdot&\cdot&\cdot\cr
13 &\cdot &\cdot &\cdot &\cdot &+&-&\cdot&\cdot&\cdot&\cr
14 &\cdot &\cdot &\cdot &\cdot &\cdot &-&+&-&+\cr
\,\,\cdot     &\cdot &\cdot &\cdot &\cdot &\cdot &\cdot &\cdot &\cdot &+\cr
\,\,\cdot     &\cdot &\cdot &\cdot &\cdot &\cdot &\cdot &\cdot &\cdot &+
} 
\]

We next show that $Q$ satisfies the evenness condition. 
With this aim, firstly it is convenient to note that for the 
considered sign pattern, in any allowed lattice path,  
the pairs of signs on the extreme positions of the horizontal 
segments alternate between $(+,-)$ and $(-,+)$.

We start by showing that for any 
$i',i''\in \{1,\ldots ,p\}\setminus I$ 
the number of indices in $I$ between $i'$ and $i''$ is even. 
To see this, assume that $\{i_s, \dots ,i_q\}$ is a maximal 
sequence of consecutive indices in $I$ between $i'$ and $i''$.   
By consecutive sequence of indices we mean 
that $i_{r+1}=i_r+1$ for $s\leq r\leq q-1$. 
Then, by Conditions~\eqref{item:C1} and~\eqref{item:C3} we must have 
$(\epsilon_{i_s j_s}\epsilon_{i_s j_{s+1}})=(+,-)$ 
because $i_s-1\not \in I$. In the same way, 
since $i_q+1 \not \in I$, 
from Conditions~\eqref{item:C2} and~\eqref{item:C3} it follows that 
$(\epsilon_{i_{q} j_{q}}\epsilon_{i_{q} j_{q+1}})=(-,+)$. 
This implies that the number of elements in $\{i_s, \dots ,i_q\}$ is even 
because the pairs of signs on the extreme positions of the horizontal segments 
of $\sP$ alternate between $(+,-)$ and $(-,+)$.  
This means that the number of elements in $Q$ between $i'+d$ and $i''+d$ 
can be expressed as a sum of even numbers, and therefore it is also even.  

Analogously, for any $j_r,j_s \in J$, 
there is an even number of elements in $\{1,\ldots ,d\}\setminus J$ 
between $j_r$ and $j_s$. 
Indeed, note that if $ \epsilon_{i_r j_r}=+1$ 
(resp.\ $\epsilon_{i_r j_r}=-1$), then by Condition~\eqref{item:C4} we have
$\epsilon_{i_r j_{r+1}}=-1$ (resp.\ $\epsilon_{i_r j_{r+1}}=+1$), 
which means that the number of elements in $\{1,\ldots ,d\}\setminus J$ 
between $j_r$ and $j_{r+1}$ is even, 
because in the considered sign pattern the signs alternate 
between $+1$ and $-1$ on each row. Therefore, 
the number of elements in $Q$ between $d-j_r+1$ and $d-j_s+1$ 
can be expressed as a sum of even numbers, and thus it is also even.  

Finally, consider the case $i=i'+d $ and 
$j=d-j_s+1$ for some $j_s\in J$ and $i'\not \in I$. 
We claim that the number of elements 
in $Q$ between $i$ and $j$ is even.  
Note that if $1\not \in I$, thanks to the previous results, 
it suffices to show that $j_1$ is odd, 
but this follows from Condition~\eqref{item:C1} 
because we must have $(-1)^{j_1-1}=+1$.  
On the other hand, if $1\in I$, let $\{i_1,\ldots ,i_s\}$ 
be the maximal sequence of consecutive indices in $I$ containing $i_1=1$. 
Due to the results above, note that to prove our claim, 
it is enough to show that $i_s+j_1-1$ is even. 
To prove this, we consider two cases. 
Assume first that $(\epsilon_{i_1 j_1},\epsilon_{i_1 j_2})=(+,-)$.  
Then, $j_1$ is odd because $(-1)^{j_1-1}=+1$. 
Since $i_s+1\not \in I$, from Conditions~\eqref{item:C2} and~\eqref{item:C3} 
it follows that the pair of signs 
$(\epsilon_{i_s j_s},\epsilon_{i_s j_{s+1}})$ 
can never be equal to $(+,-)$. Therefore, 
we conclude that $i_s$ is even,  
which means that $i_s+j_1-1 $ is also even. 
Assume now that $(\epsilon_{i_1 j_1},\epsilon_{i_1 j_2})=(-,+)$. Then, 
$j_1$ is even and, like in the previous case, 
by Conditions~\eqref{item:C2} and~\eqref{item:C3} the pair of signs 
$(\epsilon_{i_s j_s},\epsilon_{i_s j_{s+1}})$ 
can never be equal to $(+,-)$. Therefore, 
we conclude that $i_s$ is odd, which means that 
$i_s+j_1-1$ is even. 

In consequence, $Q$ satisfies Gale's evenness condition. 

Conversely, let $Q$ be a subset of $\{1,\ldots ,p+d \}$ 
with $d-1$ elements which satisfies Gale's evenness condition. 
Define the sets $I:=\{i-d \mid i\in Q , i>d\}$ 
and $J:=\{d-j+1 \mid j\not \in Q, j\leq d \}$. 
Since $Q$ has $d-1$ elements, we can write 
$I=\left\{i_1,\ldots,i_k\right\}$ and 
$J=\left\{ j_1,\dots ,j_{k+1}\right\}$ for some $k\leq d-1$, 
where we assume that $i_1<\cdots <i_k$ and $j_1<\cdots <j_{k+1}$. 
Then, with $Q$ we associate the lattice path $\sP$ defined by~\eqref{Path}. 
We next show that $\sP$ is allowed. We divide the proof into two cases.  

Firstly, it is convenient to note that applying the evenness 
condition to $d-j_r+1$ and $d-j_{r+1}+1$, for $r\leq k$, 
it follows that the columns $j_r$ and $j_{r+1}$ of $(\epsilon_{ij})$ 
always have opposite signs. 

Assume first that $1\not \in I$. Then, 
by considering $i=d+1$ and $j=d-j_1+1$ in the evenness condition, 
we conclude that $j_1$ must be odd. If $k\geq 1$, 
since the columns $j_1$ and $j_2$ have opposite signs, 
we know that $j_2$ must be even, and thus 
$(\epsilon_{i_1 j_1},\epsilon_{i_1 j_2})=(+,-)$.
Let $\{i_1,\ldots ,i_s\}$ be the maximal sequence of consecutive 
indices in $I$ containing $i_1$. Then, 
using the fact the columns $j_r$ and $j_{r+1}$ 
always have opposite signs, it follows that for $r\leq s$, 
$(\epsilon_{i_r j_r},\epsilon_{i_r j_{r+1}})=(+,-)$ if $r$ is 
odd and $(\epsilon_{i_r j_r},\epsilon_{i_r j_{r+1}})=(-,+)$ if $r$ is even.
If $i_s<p$, by the evenness condition there is an even number of elements in 
$Q$ between $i_1+d-1$ and $i_s+d+1$, so $i_s-i_1$ must be odd.  
Therefore, we have $(\epsilon_{i_s j_s},\epsilon_{i_s j_{s+1}})=(-,+)$. 
If $s\neq k$, we can repeat this argument considering $i_{s+1}$ 
instead of $i_1$. Then, after a finite number of steps, 
we either conclude that $i_k=p$ 
or that $(\epsilon_{i_k j_k},\epsilon_{i_k j_{k+1}})=(-,+)$. 
This proves that $\sP$ is allowed. 

Now assume that $1\in I$. Let $\{i_1,\ldots ,i_s\}$ be the maximal 
sequence of consecutive indices in $I$ containing $i_1=1$.  
Considering the evenness condition with $i=d+i_s+1$ and $j=d-j_1+1$, 
we conclude that $i_s+j_1$ must be odd. If $j_1$ is odd and $i_s$ is even, 
for $r\leq s$ we have $(\epsilon_{i_r j_r},\epsilon_{i_r j_{r+1}})=(+,-)$ 
for $r$ odd and $(\epsilon_{i_r j_r},\epsilon_{i_r j_{r+1}})=(-,+)$ 
for $r$ even, so in particular 
$(\epsilon_{i_s j_s},\epsilon_{i_s j_{s+1}})=(-,+)$. 
If $j_1$ is even and $i_s$ is odd, it follows that for $r\leq s$, 
$(\epsilon_{i_r j_r},\epsilon_{i_r j_{r+1}})=(-,+)$ if $r$ is 
odd and $(\epsilon_{i_r j_r},\epsilon_{i_r j_{r+1}})=(+,-)$ if $r$ is even. 
Therefore, again we have 
$(\epsilon_{i_s j_s},\epsilon_{i_s j_{s+1}})=(-,+)$. 
In both cases, we can now apply the argument used in the last 
part of the case $1\not \in I$ to show that $\sP$ is allowed. 
This concludes the proof of the theorem.
\end{proof}

\section{Tropical Half-spaces in General Position}

Recall that a $k\times k$ matrix $M$ with entries in $\rmax$ is tropically
non-singular if the tropical permanent
\[
\operatorname{tper}M:=\max_{\sigma\in S_k} \sum_{1\leq i\leq k} M_{i\sigma(i)}
\]
is finite and is attained by exactly one permutation $\sigma$. 

Consider a tropical polyhedral cone defined by the system of
inequalities $Ax\leq Bx$, where $A,B$ are $p\times d$ matrices
with entries in $\rmax$ which we may require to satisfy
$A_{ij}B_{ij}=\zero$ for all $1\leq i\leq p$ and $1\leq j\leq d$.
Let $C:=A\oplus B$. We say that the latter inequalities are {\em in general position} if any $k\times k$ submatrix of $C$ is tropically non-singular.

An elementary part of the proof of McMullen's upper bound theorem
is to show that 
the number of facets of a polytope with $p$ vertices in dimension $d$ 
is maximized by a simplicial polytope. The following theorem 
may be thought of as a tropical version of the dual of this result.

\begin{theorem}\label{th-general}
The maximal number of extreme rays of a tropical cone defined as the intersection of $p$ tropical half-spaces in dimension $d$ is attained when these half-spaces are in general position.
\end{theorem}
\begin{proof}
The proof is similar in its spirit to the one of Theorem~\ref{th-upperbound}.
We can choose a sequence of perturbed matrices 
\begin{align}\label{e-inclusion}
A(m)\leq A\qquad \text{and}\qquad B(m)\geq B \enspace
,\qquad m\in \N \enspace ,
\end{align}
in such a way that for all $m$, every square submatrix of the matrix 
$C(m):=A(m)\oplus B(m)$ is tropically non-singular and $A(m)\to A$, 
$B(m)\to B$ as $m$ tends to infinity. For instance, if $B_{ij}>-\infty$, 
we may require that $(B(m))_{ij}>B_{ij}$ and $(A(m))_{ij}=-\infty$, 
whereas if $A_{ij}>-\infty$, we may require that $(A(m))_{ij}<A_{ij}$ 
and $(B(m))_{ij}=-\infty$. The matrices 
$A(m)$ and $B(m)$ may be chosen arbitrarily close to $A$ and $B$, respectively, and if their entries are rationally independent, every submatrix of $C(m)$ must be tropically non-singular. Let $\mathcal{C}:=\{x\in \rmax^d\mid Ax\leq Bx\}$ and $\mathcal{C}(m):=\{x\in \rmax^d\mid A(m)x\leq B(m)x\}$. Due to Property~\eqref{e-inclusion}, we have $\mathcal{C}(m)\supset \mathcal{C}$.

Let $K(p,d)$ denote the maximal number of extreme rays of (tropical) polyhedral cones in dimension $d$ defined by systems of $p$ inequalities in general position. For all $m\in \N$,
let $\{u_k(m)\}_{k=1,\dots ,K(m)}$ denote a generating family of $\mathcal{C}(m)$,
which by Theorem~\ref{TheoMinkowski} can be obtained by selecting precisely one element in each extreme
ray of $\mathcal{C}(m)$, so that $K(m)\leq K(p,d)$. Possibly after extracting
a subsequence, we may assume that $K:=K(m)$ is independent of $m$. Every vector
$u_k(m)$ can be chosen to be normalized (e.g. to have the maximum of its entries
equal to $\unit$) and so, perhaps after extracting again a subsequence, we may
assume that $u_k(m)$ has a limit $u_k\in \rmax^d$ different from the zero vector $\zero$ as
$m$ tends to infinity. Since $\mathcal{C}\subset \mathcal{C}(m)$, we
deduce that for all vectors $v\in \mathcal{C}$,
and for all $m\in \N$, we can find some scalars $\lambda_k(m)$ 
such that 
\[
v=\bigoplus_{1\leq k\leq K} \lambda_k(m) u_k(m) \enspace .
\]
Since every $u_k(m)$ has some entry $i$ (depending on $k$ and $m$)
equal to $\unit$, we deduce that $\lambda_k(m)\leq v_i$, and
so $\lambda_k(m)\leq \max_j v_j$ for all $k$ and $m$. Hence,
$\lambda_k(m)$, which is bounded as $m$ tends to infinity, must have
an accumulation point $\lambda_k\in \rmax$, and we deduce that
\[
v=\bigoplus_{1\leq k\leq K}\lambda_k u_k \enspace .
\]
Moreover, by passing to the limit in $A(m)u_k(m)\leq B(m)u_k(m)$,
we deduce that $Au_k\leq Bu_k$, showing that $u_k\in \mathcal{C}$.
It follows that $\{u_k\}_{1\leq k\leq K}$ is a generating family of $\mathcal{C}$. Since the number of extreme rays of a polyhedral cone is bounded by the cardinality of any of its generating families, we deduce that the number of extreme
rays of $\mathcal{C}$ is bounded by $K(p,d)$.
\end{proof}

We denote by $\ntrop(\epsilon)$ (resp.\ $\nclass(\epsilon)$)
the number of tropically (resp.\ non-tropically) allowed
lattice paths for the sign pattern $\epsilon$.
We also denote by $\ntropext(p,d)$ the maximal number of extreme
rays of a tropical cone in dimension $d$ defined as the intersection
of $p$ half-spaces. We have shown that
\begin{align}\label{e-synthesis}
\max_{\epsilon\in \{\unit,\ominus \unit\}^{p\times d} }
\ntrop(\epsilon) \leq \ntropext(p,d)
\leq U(p+d,d-1) = 
\max_{\epsilon\in \{\pm 1\}^{p\times d}} \nclass(\epsilon) \enspace .
\end{align}
The following conjecture, which is suggested by the analogy with the classical
case, states that the two leftmost terms in the latter expression are equal.

\begin{conj}\label{conj-agk}
The number of extreme rays of a tropical cone in dimension $d$ defined 
as the intersection of $p$ half-spaces is maximized by the polar of a 
signed cyclic polyhedral cone. 
\end{conj}
This conjecture is also (weakly) supported by Theorem~\ref{th-general}:
the signed cyclic polyhedral cones can be easily seen to define systems
of constraints in general positions, and among these, they somehow
provide the simplest model.

If it were true, this conjecture would have surprising consequences
in terms of complexity, showing that polyhedra defined by $p$ constraints
in a fixed dimension $d$ have fewer extreme points in
the tropical case as $p$ tends to infinity,
see Remark~\ref{rk-paradoxal} below.

\section{Computing the number of tropically allowed paths}
\newcommand{\mdown}{\mathsf{d}}
\newcommand{\mright}{\mathsf{r}}
We next give an inductive formula allowing
one to compute the number $\ntrop(\epsilon)$ of tropically allowed lattice
paths for the sign pattern $\epsilon=(\epsilon_{ij})$ in a time
which is linear in the size of the pattern.

First, we write the signs $\epsilon_{ij}$, $1\leq i\leq p, 1\leq j\leq d$ 
in a $p\times d$ table, that we complete by adding one 
dummy row at the top numbered $0$ and one dummy row
at the bottom numbered $p+1$. 

We shall consider paths starting from
the position $(0,1)$ (row $0$, column $1$)
and ending at some position $(p+1,j)$ (row $p+1$, column $j$). 
Such paths are said to be tropically allowed if the
subpath lying in rows $1,\ldots,p$ is tropically allowed.

We represent every lattice path by a word in the alphabet $\{\mdown,\mright\}$.
The letter $\mdown$ represents a downward move, whereas the letter
$\mright$ represents a move to the right. (The letter $\mdown$ should not be confused
with the symbol $d$ for the dimension.)
For instance, if $p=1$ and $d=2$, 
the word $\mdown\mright\mdown$ corresponds to the path 
\[ 
\begin{array}{ccc}(0,1)\\
\mdown\downarrow \\
(1,1)&\stackrel{\displaystyle \mright}{\to}&
(1,2)\\
&&\mdown{\downarrow}\\
&&(2,2)
\end{array}
\]
\begin{figure}
\begin{center}
\begin{tikzpicture}[shorten >=1pt,>=stealth',bend angle=45,state/.style={circle,draw=black,very thick,minimum size=5ex,inner sep=0ex},transition/.style={}]
\node[state] (m) at (0,0) {$-$};
\node[state,accepting] (mp) at (4.5,0) {$-+$};
\node[state] (p) at (0,2.5) {$+$};
\node[state,accepting] (pm) at (4.5,2.5) {$+-$};
\node[initial above,state] (one) at (2.25,4.5) {$1$};

\draw[transition] (one) edge[loop left] node {$\mathsf{r}$} (one);
\draw[transition] (p) edge[loop left] node {$\mathsf{r}$} (p);
\draw[transition] (pm) edge[loop right] node {$+\mathsf{d}$} (pm);
\draw[transition] (m) edge[loop left] node {$\mathsf{r}$} (m);
\draw[transition] (mp) edge[loop right] node {$+\mathsf{d}$} (mp);
\draw[transition,bend angle=30] (one) edge[->,bend left] node[above right] {$\mathsf{d}$} (pm);
\draw[transition,bend angle=25] (p) edge[->,bend left] node[below] {$-\mathsf{d}$} (pm)
  (pm) edge[->,bend left] node[above] {$+\mathsf{r}$} (p);
\draw[transition] (pm.south west) edge[->] node[right=1ex] {$-\mathsf{r}$} (m.north east);
\draw[transition,bend angle=25] (m) edge[->,bend left] node[below] {$+\mathsf{d}$} (mp) 
  (mp) edge[->,bend left] node[above] {$-\mathsf{r}$} (m);
\end{tikzpicture}
\end{center}
\caption{An automaton recognizing tropically allowed paths.} 
\label{FigAutomaton}
\end{figure}
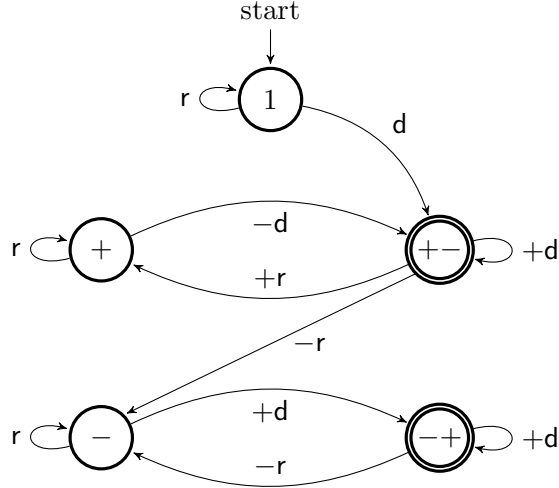
Consider now the automaton represented in Figure~\ref{FigAutomaton},
in which the state denoted by $1$ (with an incoming arrow) is initial
and the states denoted by the symbols $+-$ and $-+$
(with double circles) are final.

The arcs are labeled
by letters, and sometimes by signs. We next introduce
an acceptance condition which slightly differs from the classical
one in automata theory, in order to take into account the sign
pattern.

A word is said to be {\em accepted} by the automaton
if the following holds. 
We read the letters of the word from left to right, performing
at the same time the corresponding moves (downward, or to the right) 
in the table and in the automaton (following arcs).
A move is accepted only if the sign $\epsilon_{ij}$ of the current position
of the table is the same as the sign of the corresponding arc originating from the current state on the automaton (if there
is no sign on this arc, $\epsilon_{ij}$ can be arbitrary). The word
is accepted if, when starting from position $(0,1)$ in the table and
from the initial node in the automaton, every successive move is accepted, leading to a state of the automaton which is final, and if the final position in the table is at some point $(p+1,j)$ with $1\leq j\leq d$, which means that the word contains precisely $p+1$ occurrence of the letter $\mdown$ and at most $d-1$ occurrences of the letter $\mright$.

For instance, if $p=1$ and $d=2$, and if the sign pattern is $[+,-]$, 
the word $\mdown\mdown$ is accepted since it corresponds
to the following path in the automaton:
\[ 1\stackrel{\mdown}{\to}+-\stackrel{+\mdown}{\to}+- 
\]
Similarly, the word $\mdown\mright\mdown$ is accepted, since
it corresponds to the path:
\[1\stackrel{\mdown}{\to}+-\stackrel{+\mright}{\to}+ \stackrel{-\mdown}{\to}+-
\]

The introduction of the previous automaton is motivated by the following
result.
\begin{prop}\label{prop-autom}
The tropically allowed lattice paths are in one to one correspondence
with the words that are accepted by the automaton, and each of these words
corresponds precisely to one path in the automaton.
\end{prop}
\begin{proof}
Imagine a pen, drawing the path starting  from the top left position and
making only moves downward or to the right.
We shall see that the states of the automaton are used to
record the information necessary to determine how the pen can be moved
to draw a tropically allowed lattice path. 

First, the pen is at position $(0,1)$ (on the dummy top row),
and the current state of the automaton
is the initial state, $1$. Then, the pen may either stay
on the dummy top row, moving to the right, or leave the dummy row,
moving down, corresponding to the two arcs $1\stackrel{\mright}{\to} 1$
and $1\stackrel{\mdown}{\to}{+-}$. Assume that the latter arc has been chosen
after a sequence of moves to the right (which cannot exceed $d-1$ due to the
final acceptance condition), so that the pen is now at some position
$(1,j)$ with $1\leq j\leq d$. Then, the pen always may move to the right,
beginning an horizontal segment. If $\epsilon_{1j}=\oplus \unit$,
in accordance with Condition~\eqref{item:C1}, the pen may also move down.
These moves correspond to the three arcs leaving node $+-$ in the automaton:
we use the state $+$ (resp.\ $-$) to record that the horizontal segment which has been opened starts with a $+$ (resp.\ $-$ sign). 

Consider now the situation in which
$\epsilon_{1j}=+1$ and a move to the right has been selected, so
that the current state in the automaton is $+$ and the position 
of the pen is now $(1,j+1)$. Since the sign of every position
of an horizontal segment which is not extreme does not matter in
the definition of tropically allowed path, the move to the right
can always be selected. By Condition~\eqref{item:C4}, a downward move
can be accepted only if the sign at the current position is $-$,
since an horizontal segment which began with a $+$ must end by a $-$,
and since the downward move ends the current horizontal segment.
The latter move corresponds to the arc $+\stackrel{-\mdown}{\to}+-$
in the automaton. 

Similarly, the state $-$ indicates
that the pen is now drawing an horizontal segment starting from a $-$ sign,
and the state $-+$ indicates that such a segment has been
closed. Observe that there is an arc from state $+-$ to state $-$, but no arc
from state $-+$ to state $+$, because, by Condition~\eqref{item:C5}, the
pair $(-,+)$ may always appear after a pair $(+,-)$ as the
signs of the extreme positions of an horizontal segment, whereas
the opposite is not allowed. 

With this interpretation in mind, it is readily seen that every accepted
word bijectively corresponds to a tropically allowed lattice path.

An inspection of the automaton
also shows that it is unambiguous, meaning that there is precisely
one path in the automaton for each accepted word. Indeed, the unambiguity
stems from the fact that at each state, there is at most one leaving
arc with a given letter and sign.
\end{proof}

The inductive formula to compute $\ntrop(\epsilon)$ is next obtained
by some elementary bookkeeping.

We denote by $\chi^+(i,j)$ the number which is $1$ if
$\epsilon_{ij}=+ 1$ and $0$ otherwise. Similarly, 
$\chi^-(i,j)$ is $1$ if $\epsilon_{ij}=- 1$ and $0$ otherwise. 
For $0\leq i\leq p+1,1\leq j\leq d+1$, define the numbers
$N_+(i,j)$, $N_-(i,j)$, $N_{+-}(i,j)$, 
$N_{-+}(i,j)$, and $N_1(0,j)$ by the following 
inductive formul\ae
\begin{align*}
N_1(0,j)& = N_{1}(0,j+1)+ N_{+-}(1,j) \enspace , && 1\leq j\leq d \enspace , \\
N_+(i,j)&= N_+(i,j+1)+\chi^-(i,j)N_{+-}(i+1,j) \enspace , && 0\leq i\leq p, 1\leq j\leq d \enspace , \\
N_-(i,j)&= N_-(i,j+1)+\chi^+(i,j)N_{-+}(i+1,j) \enspace , && 0\leq i\leq p, 1\leq j\leq d \enspace ,\\
N_{+-}(i,j)&= \chi^+(i,j)N_{+-}(i+1,j)+\chi^+(i,j)N_{+}(i,j+1) \\
& \phantom{={}} \qquad {}+\chi^-(i,j)N_-(i,j+1)\enspace , && 0\leq i\leq p, 1\leq j\leq d \enspace , \\
N_{-+}(i,j)&= \chi^-(i,j)N_{-}(i,j+1)+\chi^+(i,j)N_{-+}(i+1,j) \enspace , && 0\leq i\leq p, 1\leq j\leq d\enspace ,
\end{align*}
together with the boundary conditions
\begin{align*}
N_{s}(i,d+1)&= 0 \enspace , && 0\leq i\leq p+1 \enspace ,\; s\in \{+,-,+-,-+\} \enspace , \\
N_{1}(0,d+1)&= 0 \enspace , \\
N_{s}(p+1,j)&= 1 \enspace , && 1\leq j\leq d \enspace ,\; s\in \{+-,-+\} \enspace , \\
N_{s}(p+1,j)&= 0 \enspace , && 1\leq j\leq d \enspace ,\; s\in \{+,-\}\enspace .\end{align*}
\begin{coro}[Computing the number of tropically allowed paths]\label{e-compute}
For all sign patterns $\epsilon$, we have
\[
\ntrop(\epsilon)=N_1(0,0) \enspace .
\]
\end{coro}
\begin{proof}
We claim that for each state $s$ of the automaton, and for all $0\leq i\leq p$, $1\leq j\leq d$, $N_{s}(i,j)$ represents the number of possible sequences of remaining moves of a pen drawing a tropically
allowed path, given that the current position of the pen is $(i,j)$ and that
the previous moves of the pen led to this position and to state $s$.

We observe that the equations above, except for the two ones which determine the boundary values $N_{s}(p+1,j)$, are readily obtained from the automaton.
For instance, the formula for $N_{+-}(i,j)$ as a sum of three terms corresponds
to the three options: move down if the sign $\epsilon_{ij}$ is positive;
open an horizontal segment with initial sign $+$ under the same
condition; or open an horizontal segment with initial sign $-$
if $\epsilon_{ij}$ is negative. The other formul\ae\ are obtained
in a similar way. Note that the boundary conditions which determine
$N_{s}(p+1,j)$ force the final state to be either $+-$
or $-+$, meaning that every horizontal path which has been opened
must have been closed. Using these considerations, one readily shows
the claim by a backward induction on $(i,j)$, initialized
when $i=p+1$ or $j=d+1$.
\end{proof}

\begin{remark}
To compute the number of (non-tropically) allowed paths, 
it suffices to add an arc $-+\stackrel{+\mright}{\to}+$ in the automaton.
Then,
we must add a third term $\chi^+(i,j)N_{+}(i,j+1)$ in the expression
of $N_{-+}(i,j)$, and one can check that the number $N_1(0,0)$ now determines
the number of allowed paths. One can also check that $N_{-+}(i,j)=N_{+-}(i,j)$,
meaning that the automaton is no longer minimal (the states $+-$ and $-+$ 
can be identified).
\end{remark}

\section{Upper and lower estimates for the number of extreme rays of the polar of signed cyclic polyhedral cones}
\label{sec-estimates}

We showed that $\ntropext(p,d)$, the maximal number of extreme rays of a tropical polyhedral cone defined by $p$ inequalities in dimension $d$ is bounded from above by its classical analogue, $U(p+d,d-1)$, and bounded from below by $\ntrop(p,d)$, the maximal number
of tropically allowed lattice paths for a $p\times d$ signed pattern, see
Eqn~\eqref{e-synthesis}. The asymptotic behavior of $U(p,d)$ is easily determined. In this section, we provide explicit estimates for the lower bound $\ntrop(p,d)$ and derive its asymptotic behavior as $p$ or $d$ tends to infinity.

We shall say that a tropically allowed path is of 
$\pmp $ type if the pair of signs consisting of the
signs of the leftmost and rightmost positions of 
each of its horizontal segments is $(-,+)$. 
Tropically allowed paths of $\ppm $ type are defined in a symmetric way. 
Recall that a tropically allowed path consists of a path of 
$\ppm $ type followed by a path of $\pmp $ type, 
one of these being possibly empty.

Let $\nmp(p,d)$ (resp.\ $\npm(p,d)$) 
denote the maximal number of tropically allowed paths of $\pmp $ type 
(resp.\ $\ppm $ type) in a $p\times d$ sign pattern. 
We shall also need $\npm_\ell(p,d)$, 
which denotes the maximal number of tropically allowed paths of 
$\ppm $ type using the last column of a $p\times d$ sign pattern.
We make the following observation:
\begin{align}
\nmp(p,d)=\npm(p,d) \enspace .
\end{align}
Indeed, if we read a tropically allowed path of $\pmp $ type for a 
$p\times d$ sign pattern in a reverse way (starting from the end), 
it becomes a tropically allowed path of $\ppm $ type in the reversed sign 
pattern (in which the bottom right corner becomes the top left corner), 
and vice versa. 

We first bound $\ntrop(p,d)$ from above.

\begin{prop}\label{prop-upperbound}
For every $p,d$,  
\begin{equation}\label{UpperBoundCyclic}
\ntrop(p,d) \leq (p (d-1)+1) 2^{d-1} \enspace .
\end{equation}
\end{prop}

\begin{proof}
In the first place, we claim that
\begin{align}\label{NTropPathUpperBound}
\ntrop(p,d) \leq \left( \sum_{r = 1}^p \sum_{m = 1}^{d-1} 
\npm_\ell(r-1,m) \nmp(p-r,d-m) \right) + \npm(p,d) \enspace .
\end{align}
Indeed, in this expression $(r,m)$ represents the 
leftmost position of the first $\pmp $ segment, if any, 
of a tropically allowed path for a given sign pattern. 
Then, the part of the path before this segment must be of $\ppm $ type 
in the $(r-1)\times m$ upper left submatrix of which it uses the last column, 
accounting for the term $\npm_\ell(r-1,m)$, 
whereas the part of the path after this segment must be of $\pmp $ 
type in the $(p-r)\times (d-m)$ bottom right submatrix. 
The term outside the parenthesis represents the paths which are purely of $\ppm $ type.  
The case $m=d$ is excluded because $(r,m)$ is 
supposed to be the leftmost position of a $\pmp $ segment, 
so it cannot belong to the last column. 

We claim that, for every $p$,
\begin{align}\label{BoundNmp}
\nmp(p,d)\leq 2^d -1 \enspace .
\end{align}
To see this, let $j_1,\ldots,j_{k+1}$ denote the columns used by a tropically
allowed path of $\pmp $ type in a $p\times d$ sign pattern $(\epsilon_{ij})$. 
This path is uniquely determined by $j_1,\ldots,j_{k+1}$ because due to  
Conditions~\eqref{item:C1} and~\eqref{item:C3}, the vertical ordinates $i_1,\ldots,i_k$ 
of its horizontal segments are given recursively by 
$i_1=\min\left\{i\mid \epsilon_{i j_1}=\ominus \unit \right\}$ and 
$i_r = \min\left\{i>i_{r-1}\mid \epsilon_{i j_{r}}=\ominus \unit \right\}$, 
for $r=2,\ldots ,k$. 
Since $\{j_1,\ldots,j_{k+1}\}$ can be any non-empty subset of 
$\{1,\ldots,d\}$, the bound~\eqref{BoundNmp} follows.
A similar argument shows that
\[
\npm_\ell(p,d)\leq 2^{d-1} \; 
\]
because in this case $d$ always belongs to $\{j_1,\ldots,j_{k+1}\}$, 
so $j_{k+1}=d$ and $\{j_1,\ldots,j_k\}$ can be any subset of 
$\{1,\ldots,{d-1}\}$. 

Collecting the previous bounds and using the fact that $\nmp(p,d)=\npm(p,d)$, 
from~\eqref{NTropPathUpperBound} we obtain
\begin{align*}
\ntrop(p,d) \leq p \sum_{1\leq m\leq d-1} 2^{m-1}(2^{d-m}-1) + 
2^d-1 \enspace , 
\end{align*}
which implies~\eqref{UpperBoundCyclic}.
\end{proof}

The following propositions provide lower bounds for the maximal number 
of tropically allowed paths in a $p\times d$ sign pattern. 

\begin{prop}\label{prop-lowerbound}
For $p\geq 2d$, we have 
\begin{equation}\label{LowerBoundCyclic}
\ntrop(p,d) \geq (p-2d+7) (2^{d-2}-2) \enspace .
\end{equation}
\end{prop}

\begin{proof}

\begin{figure}
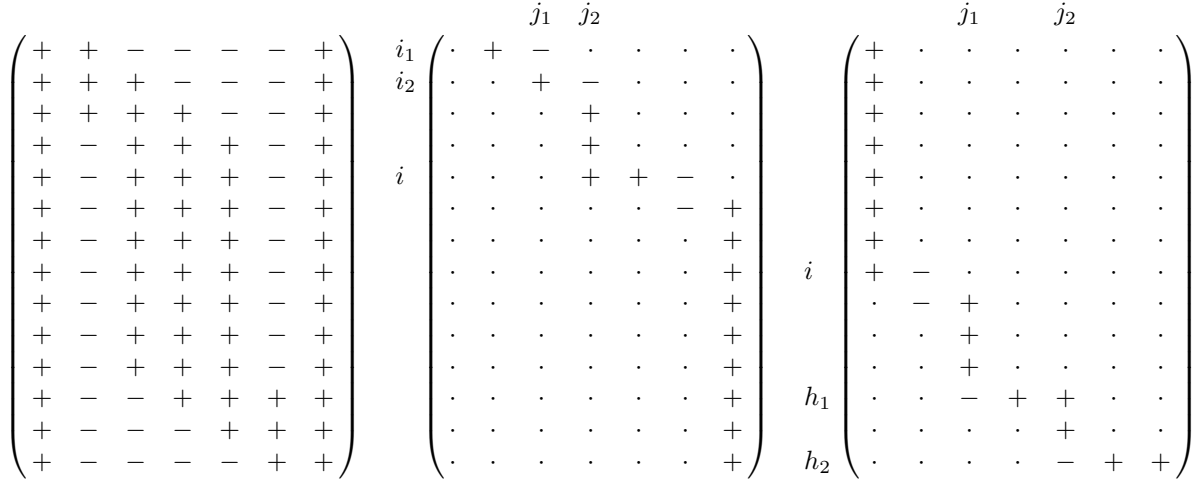

{\small 
\[
\bordermatrix{
&      &      &      &      &      &      &      \cr
& \mop & \mop & \mon & \mon & \mon & \mon & \mop \cr
& \mop & \mop & \mop & \mon & \mon & \mon & \mop \cr
& \mop & \mop & \mop & \mop & \mon & \mon & \mop \cr
& \mop & \mon & \mop & \mop & \mop & \mon & \mop \cr
& \mop & \mon & \mop & \mop & \mop & \mon & \mop \cr
& \mop & \mon & \mop & \mop & \mop & \mon & \mop \cr
& \mop & \mon & \mop & \mop & \mop & \mon & \mop \cr
& \mop & \mon & \mop & \mop & \mop & \mon & \mop \cr
& \mop & \mon & \mop & \mop & \mop & \mon & \mop \cr
& \mop & \mon & \mop & \mop & \mop & \mon & \mop \cr
& \mop & \mon & \mop & \mop & \mop & \mon & \mop \cr
& \mop & \mon & \mon & \mop & \mop & \mop & \mop \cr
& \mop & \mon & \mon & \mon & \mop & \mop & \mop \cr
& \mop & \mon & \mon & \mon & \mon & \mop & \mop 
}
\quad
\bordermatrix{
    &       &       & j_1   & j_2   &       &       &       \cr
i_1 & \cdot & \mop  & \mon  & \cdot & \cdot & \cdot & \cdot \cr
i_2 & \cdot & \cdot & \mop  & \mon  & \cdot & \cdot & \cdot \cr
    & \cdot & \cdot & \cdot & \mop  & \cdot & \cdot & \cdot \cr
    & \cdot & \cdot & \cdot & \mop  & \cdot & \cdot & \cdot \cr
i   & \cdot & \cdot & \cdot & \mop  & \mop  & \mon  & \cdot \cr
    & \cdot & \cdot & \cdot & \cdot & \cdot & \mon  & \mop  \cr
    & \cdot & \cdot & \cdot & \cdot & \cdot & \cdot & \mop  \cr
    & \cdot & \cdot & \cdot & \cdot & \cdot & \cdot & \mop  \cr
    & \cdot & \cdot & \cdot & \cdot & \cdot & \cdot & \mop  \cr
    & \cdot & \cdot & \cdot & \cdot & \cdot & \cdot & \mop  \cr
    & \cdot & \cdot & \cdot & \cdot & \cdot & \cdot & \mop  \cr
    & \cdot & \cdot & \cdot & \cdot & \cdot & \cdot & \mop  \cr
    & \cdot & \cdot & \cdot & \cdot & \cdot & \cdot & \mop  \cr
    & \cdot & \cdot & \cdot & \cdot & \cdot & \cdot & \mop 
}
\quad
\bordermatrix{
    &       &       & j_1   &       & j_2   &       &       \cr
    & \mop  & \cdot & \cdot & \cdot & \cdot & \cdot & \cdot \cr
    & \mop  & \cdot & \cdot & \cdot & \cdot & \cdot & \cdot \cr
    & \mop  & \cdot & \cdot & \cdot & \cdot & \cdot & \cdot \cr
    & \mop  & \cdot & \cdot & \cdot & \cdot & \cdot & \cdot \cr
    & \mop  & \cdot & \cdot & \cdot & \cdot & \cdot & \cdot \cr
    & \mop  & \cdot & \cdot & \cdot & \cdot & \cdot & \cdot \cr
    & \mop  & \cdot & \cdot & \cdot & \cdot & \cdot & \cdot \cr
i   & \mop  & \mon  & \cdot & \cdot & \cdot & \cdot & \cdot \cr
    & \cdot & \mon  & \mop  & \cdot & \cdot & \cdot & \cdot \cr
    & \cdot & \cdot & \mop  & \cdot & \cdot & \cdot & \cdot \cr
    & \cdot & \cdot & \mop  & \cdot & \cdot & \cdot & \cdot \cr
h_1 & \cdot & \cdot & \mon  & \mop  & \mop  & \cdot & \cdot \cr
    & \cdot & \cdot & \cdot & \cdot & \mop  & \cdot & \cdot \cr
h_2 & \cdot & \cdot & \cdot & \cdot & \mon  & \mop  & \mop 
}
\]
}
\caption{Sign pattern with a natural symbol shape and two tropically allowed paths.}\label{FigBecarre}
\end{figure}

We shall give a $p\times d$ sign pattern which has at least 
$(p-2d+7)(2^{d-2}-2)$ tropically allowed paths.  
Consider the $p\times d$ sign pattern $(\epsilon_{ij})$,
with a natural symbol shape ($\natural$), defined as follows: 
\[
\epsilon_{ij}=\ominus \unit \iff
\left\{
\begin{array}{l}
i=2 \text{ and } j\geq d-3 \enspace , \\
i=d-1 \text{ and } j\leq p-d+4 \enspace , \\ 
3\leq i \leq d-2 \text{ and } j\leq i-2 \enspace , \\ 
3\leq i \leq d-2 \text{ and } j\geq i+p-d+2 \enspace . 
\end{array}
\right. 
\]
An example for $p=14$ and $d=7$ is given on the left hand side of 
Figure~\ref{FigBecarre}.   
 
Let $\left\{ j_1,\ldots ,j_k \right\}$ be any non-empty subset of 
$\left\{ 3,\ldots ,d-2 \right\} $ and $i\in \left\{ d-3,\ldots ,p-d+3 \right\}$. 
Then, it can be checked that the following lattice paths
\begin{align*}
& (1,2),\ldots ,(i_1,2),\ldots ,(i_1,j_1),\ldots ,(i_k,j_k),\ldots ,(i,j_k),\ldots ,(i,d-1),(i+1,d-1),(i+1,d),\dots ,(p,d) \\
& (1,1),\dots ,(i,1),(i,2),(i+1,2),\dots ,(i+1,j_1),\dots ,(h_1,j_1),\dots ,(h_k,j_k),\dots,\!(h_k,d-1),\dots,\!(p,d-1) \\
& (1,1),\ldots ,(i_1,1),\ldots ,(i_1,j_1),\ldots ,(i_k,j_k),\ldots ,(i,j_k),\ldots ,(i,d-1),(i+1,d-1),(i+1,d),\dots ,(p,d) \\
& (1,1),\ldots ,(i,1),(i,2),(i+1,2),\ldots ,(i+1,j_1),\ldots ,(h_1,j_1),\ldots ,(h_k,j_k),\ldots ,(h_k,d),\ldots ,(p,d)
\end{align*}
where $i_r=j_r-2$ and $h_r=j_r+p-d+2$ for $r=1,\ldots ,k$, are tropically allowed. 
Examples of the first and last cases are given in Figure~\ref{FigBecarre} for $k=2$. Indeed, 
note that in the last two cases $\left\{ j_1,\ldots ,j_k \right\}$ can also be empty, 
in which case these paths reduce to 
\begin{align*}
& (1,1),\ldots ,(i,1),\ldots ,(i,d-1),(i+1,d-1),(i+1,d),\dots ,(p,d) \enspace \text{ and } \\
& (1,1),\ldots ,(i,1),(i,2),(i+1,2),\ldots ,(i+1,d),\ldots ,(p,d)
\end{align*} 
respectively. Therefore, since all these paths are different, 
for this sign pattern we have at least 
$2(p-2d+7)(2^{d-4}-1)+2(p-2d+7)2^{d-4}=(p-2d+7)(2^{d-2}-2)$ 
tropically allowed paths. 
\end{proof}

\begin{prop}\label{prop-lowerboundPfixed}
For $d\geq 2p+1$, we have 
\begin{equation}\label{LowerBoundCyclicPfixed}
\ntrop(p,d) \geq U(d,d-p-1) \enspace .
\end{equation}
\end{prop}

\begin{proof}
Consider the $p\times d$ sign pattern $(\epsilon_{ij})$ defined by 
$\epsilon_{ij}:=\ominus \unit$ if and only if $i+j$ is odd. 
We shall show that for this sign pattern, when $d\geq 2p+1$, 
there exist at least $U(d,d-p-1)$ tropically allowed lattice paths. 
 
Let $Q$ be any subset of $\{1,\ldots ,d\}$ with $d-p-1$ elements 
which satisfies Gale's evenness condition, 
i.e. such that for any $j',j''\in \{1,\ldots d\}\setminus Q$ 
the number of elements in $Q$ between $j'$ and $j''$ is even. 
Assume that $\{1,\ldots d\}\setminus Q=\{j_1,\ldots , j_{p+1}\}$, 
where $j_1<\cdots <j_{p+1}$. Then, 
the lattice path 
\[
(1,j_1),(1,j_2),(2,j_2),(2,j_3),\ldots ,(p,j_p),(p,j_{p+1}) 
\]
is tropically allowed. Indeed, by Gale's evenness condition applied to 
$j'=j_r$ and $j''=j_{r+1}$, 
the signs in the positions $(r,j_r)$ and $(r,j_{r+1})$ must be opposite. 
Since the signs in the positions $(r,j_{r+1})$ and $(r+1, j_{r+1})$ 
are also opposite, we conclude that 
$(\epsilon_{rj_r},\epsilon_{rj_{r+1}})=(+,-)$ or 
$(\epsilon_{rj_r},\epsilon_{rj_{r+1}})=(-,+)$ for all $1\leq r \leq p$, 
depending on whether $j_1$ is odd or not. Therefore,  
the path above is tropically allowed. 

Since there are $U(d,d-p-1)$ subsets of $\{1,\ldots ,d\}$ with $d-p-1$ 
elements which satisfy the evenness condition, the proposition follows. 
\end{proof}

The following proposition points out cases
in which the upper bound is attained.
\begin{prop}\label{prop-attained}
The upper bound $U(p+d,d-1)$ for $\ntropext(p,d)$ is attained for $p\leq 3$,
for $d\leq 4$, and for $p=4$ and $d$ even.
\end{prop}

\begin{proof}
We shall only give the sign patterns for which the polar of the signed cyclic 
polyhedral cone attains the bound, leaving the details to the reader.

For $d\leq 4$ it is enough to 
define $\epsilon_{ij}=\ominus \unit$ if and only if $j=2$.

When $p=1$ the maximizing sign pattern is given by 
$\epsilon_{1j}=\ominus \unit$ if and only if $j$ is even. 

For $p=2$ we have to define $\epsilon_{ij}=\ominus \unit$ 
if and only if $i+j$ is odd,  
but $\epsilon_{pd}=\oplus \unit$ when $d$ is odd even if $p+d$ is odd. 

The case $p=3$ needs to be divided. 
If $d$ is even the maximizing sign pattern is given by 
$\epsilon_{ij}=\ominus \unit$ if and only if $i+j$ is odd, 
but $\epsilon_{pd}=\oplus \unit$ even if $p+d$ is odd. 
When $d$ is odd, $\epsilon_{ij}=\ominus \unit$ if and only if 
$i+j$ is even, but $\epsilon_{11}=\oplus \unit$ and 
$\epsilon_{pd}=\oplus \unit$ even if $p+d$ is even.

Finally, when $p=4$ and $d$ is even, the maximizing sign pattern 
is given by $\epsilon_{ij}=\ominus \unit$ if and only if $i+j$ is even, 
except for $\epsilon_{11}$ and $\epsilon_{pd}$ which must be equal to 
$\oplus \unit$. 
\end{proof}

\begin{remark}
The bound $U(p+d,d-1)$ can be written as
\[
\mychoose{p+%
k
}{k-1
}
+ 
\mychoose{p+%
k-1
}{
k-2
}
\qquad\text{when $d=2k-1$, and} 
\]
\[
2\mychoose{p+k
}{k-1
}
\qquad\text{when }d=2k. 
\]
\end{remark}
\begin{remark}\label{rk-paradoxal}
An interesting situation arises when the dimension $d$ is kept fixed,
whereas the number of constraints $p$ tends to infinity. Then, it follows
readily from the previous formula that
\[ U(p+d,d-1) = \Theta(p^{\floor{\frac{d-1}{2}}}) \qquad \text{as }\qquad
p\to\infty 
\]
whereas it follows from Propositions~\ref{prop-upperbound} 
and~\ref{prop-lowerbound} that 
\[
\ntrop(p,d) = \Theta(p) \qquad \text{as }\qquad
p\to\infty \enspace 
\]
(these asymptotic expansions of course are not uniform in $d$).
Hence, if Conjecture~\ref{conj-agk} was true, when the dimension $d$ is fixed,
and assuming that $d\geq 5$, the number of extreme points of a polyhedral cone
defined by $p$ constraints in dimension $d$ would grow much more slowly
in the tropical case, showing only a linear growth.
\end{remark}
\begin{remark}
When the number of constraints $p$ is kept fixed, whereas
$d$ tends to infinity,
it is easily seen that the upper bound
$U(p+d,d-1)$ for the number of extreme rays $\ntropext(p,d)$ is equivalent
to the lower bound $U(d,d-p-1)$ of Proposition~\ref{prop-lowerboundPfixed}.
It follows that
\[
\ntropext(p,d)\sim U(p+d,d-1) \qquad \text{as } d\to\infty \enspace .
\]
In other words, the inequalities in~\eqref{e-synthesis} are asymptotically
tight when $d\to\infty$.
\end{remark}

\def\ric{\text{r}}
\def\xav{\text{x}}
\def\st{\text{s}}

We illustrate the previous results by displaying, in Table~\ref{lowervsupper}, for each value of $(p,d)$ the best known bounds for $\ntropext(p,d)$.
Each entry
of the table is an interval containing $\ntropext(p,d)$.
When the upper and lower bounds coincide,
we write a number instead of the interval reduced to this number. 
The upper bounds 
come from Theorem~\ref{th-upperbound}. To get lower bounds,
we use Theorem~\ref{th-allowed}, which implies that
$\ntropext(p,d)\geq\ntrop(\epsilon)$ for all sign patterns.
Then, we consider explicit sign patterns $\epsilon$, which come
either from Proposition~\ref{prop-attained},
or from computer experiments.
Indeed, for all the values of $p,d$ such that $pd\leq 30$, we computed
$\ntrop(\epsilon)$ for the $2^{pd}$ sign patterns $\epsilon$, so that the lower bound actually gives $\ntrop(p,d)$, which is the conjectured value for $\ntropext(p,d)$. From these ``low dimensional'' cases, we derived some
plausible values for the patterns $\epsilon$ maximizing or approaching
$\ntrop(\epsilon)$  for higher values of $(p,d)$, in particular variations
on the ``natural'' pattern introduced in the proof of Proposition~\ref{prop-lowerbound}. Experiments actually indicate that there is no simple universal
maximizing sign pattern. Finding the optimal patterns (and so, computing
$\ntrop(p,d)$) seems to be an interesting combinatorial problem, 
which is beyond the scope of the present paper.

\begin{table}
{\tiny\[\def\<#1/#2(#3){[#2,#1]}\def\=<#1/#2(#3){#1}\def\?<#1?{[?,#1]}
\begin{array}{c|cccccccccccccccccc}
d\;  \backslash \;  p & 1 & 2 & 3 & 4 & 5 & 6 & 7 & 8 & 9 & 10 & 11 
\cr
\hline
3& 4& 5& 6& 7& 8& 9& 10& 11& 12& 13& 14 
\cr
4&  6& 8& 10& 12& 14& 16& 18& 20& 22& 24& 26 
%
\cr
5 & 9& 14& 20& \<27/26(\st,\xav)& \<35/32(\xav,\st)& \<44/38(\st)& \<54/44(\st)& \<65/50(\st)& \<77/56(\st)& \<90/62(\st)& \<104/68(\st) 
\cr
6& 12&    
 20& 30& \=<42/42(\ric,\xav)& \<56/55(\xav,\st)& \<72/68(\st)& \<90/82(\st)& \<110/96(\st)& \<132/110(\st)& \<156/124(\st)& \<182/138(\st) 
\cr
%
7& 16& 30& 50& \<77/71(\ric,\st)& \<112/96(\st)& \<156/124(\st)& \<210/152(\st)& \<275/180(\st-176\ric)& \<352/208(\st-204\ric)& \<442/236(\st-232\ric)& \<546/264(\st-260\ric) 
\cr
8& 20& 40& 70& \=<112/112(\ric)& \<168/159(\st)& \<240/216(\st)& \<330/280(\st)& \<440/340(\st-336ric)& \<572/401(\st-392\ric)& \<728/452(\st-448\ric)& \<910/508(\st-504\ric) 
\cr
9& 25& 55& 105& \<182/172(\st)& \<294/250(\ric)& \<450/321(\ric)& \<660/436(\ric)& \<935/613(\ric)& \<1287/751(\ric)& \<1729/869(\ric)& \<2275/981(\ric) 
\cr
10& 30& 70& 140& \=<252/252(\ric)& \<420/370(\ric)& \<660/538(\ric)& \<990/668(\ric)& \<1430/898(\ric)& \<2002/1320(\ric)& \<2730/1642(\ric)& \<3640/1902(\ric)
\cr
11& 36& 91& 196& \<378/363(\st-322\ric)& \<672/584(\ric)& \<1122/805(\ric)& \<1782/1122(\ric)& \<2717/1357(\ric)& \<4004/1799(\ric)& \<5733/2771(\ric)& \<8008/3528(\ric)
\end{array}
\]}
\caption{Lower and upper bounds for $\ntropext(p,d)$, the maximal number of extreme rays of a tropical polyhedral cone defined by $p$ inequalities in dimension $d$.}\label{lowervsupper}
\end{table}
\def\cprime{$'$}

\end{document}